\documentclass[a4paper,12pt]{amsart}

\usepackage{graphicx}
\usepackage{subfigure}
\usepackage{amssymb}
\usepackage[usenames, dvipsnames]{color}
\newcommand{\RR}{{\bf R}}

\newcommand{\ZZ}{{\bf Z}}
\newcommand{\Ztwo}{\ZZ_2}
\newcommand{\plano}{{\RR}^2}

\newcommand{\xy}{(x,y)}

\newcommand{\seta}{\longrightarrow}

\newcommand{\contr}{{\mathcal C}}

\newcommand{\GG}{{\mathcal G}}
\newcommand{\LL}{{\mathcal L}}
\newcommand{\MM}{{\mathcal M}}
\newcommand{\OO}{{\mathcal O}}
\newcommand{\PP}{{\mathcal P}}

\newcommand{\projective}{\mathbf{RP}^1}
\newcommand{\va}{\varphi}

\DeclareMathOperator{\tr}{Tr}
\DeclareMathOperator{\cod}{cod}

\newcounter{lixo}
\newcounter{cubico}

\newtheorem{theorem}{Theorem}[section]
\newtheorem{lemma}[theorem]{Lemma}
\newtheorem{proposition}[theorem]{Proposition}
\newtheorem{corollary}[theorem]{Corollary}
\newtheorem{definition}[theorem]{Definition}
\theoremstyle{definition}
\newtheorem{Ex}[theorem]{Example}
\begin{document}

\title[Dynamics with a star node and contracting nonlinearity: \today]{Global planar dynamics with a star node and contracting nonlinearity\\ \today}

\author[B.Alarc\'on]{Bego\~na Alarc\'on}
\address{B.Alarc\'on\\Departamento de Matem\'atica Aplicada, Instituto de Matem\'atica e Estat\'\i stica, Universidade Federal Fluminense \\
Rua Professor Marcos Waldemar de Freitas Reis, S/N, Bloco H \\
Campus do Gragoat\'a,
CEP 24.210 -- 201, S\~ao Domingos --- Niter\'oi, RJ, Brasil}
\email{balarcon@id.uff.br}
\thanks{The first author was partially supported by the Spanish Research Project MINECO-18-MTM2017-87697-P}
\author[S.B.S.D. Castro]{Sofia B.S.D. Castro}
\address{S.B.S.D. Castro\\Centro de Matem\'atica da Universidade do Porto\\ Rua do Campo Alegre 687, 4169-007 Porto, Portugal
and
Faculdade de Economia do Porto \\ Rua Dr. Roberto Frias, 4200-464 Porto, Portugal}
\email{sdcastro@fep.up.pt}
\thanks{The last two authors were partially supported by CMUP, member of LASI, which is financed by national  funds through  FCT --- Funda\c c\~ao para a Ci\^encia e a Tecnologia, I.P. (Portugal)  under the projects with reference  UIDB/00144/2020 and UIDP/00144/2020.}
\author[I.S. Labouriau]{Isabel S. Labouriau}
\address{I.S. Labouriau\\Centro de Matem\'atica da Universidade do Porto\\ Rua do Campo Alegre 687, 4169-007 Porto, Portugal}
\email{islabour@fc.up.pt}

\begin{abstract}
This is a complete study of the dynamics of polynomial planar vector fields whose linear part is a multiple of the identity and whose nonlinear part is a contracting homogeneous polynomial.
The contracting nonlinearity provides the existence of an invariant circle and allows us to obtain a classification through a complete invariant for the dynamics, extending previous work by other authors that was mainly concerned with the existence and number of limit cycles.
The general results are also applied to some classes of examples: definite nonlinearities, $\ZZ_2\oplus\ZZ_2$ symmetric systems and  nonlinearities of degree 3,
for which we provide  complete sets of phase portraits.
\end{abstract}

\maketitle


%
%
%{\be Coment\'arios e altera\c c\~oes da Bego\~na, desta cor}
%
%{\sof Coment\'arios e altera\c c\~oes da Sofia, desta cor}
%
%{\isl Coment\'arios e altera\c c\~oes da Isabel, desta cor}
%
%{\sai  Desta cor \'e para sair}

\section{Introduction}

Global planar dynamics of polynomial vector fields has been of interest for many years. Part of this interest arises from its connection to Hilbert's $16^{th}$ problem on the number of limit cycles for the dynamics.
Because of Hilbert's $16^{th}$ problem, substantial effort has been devoted in establishing a bound for the number of limit cycles. For some contributions in this direction when the vector field has homogeneous nonlinearities see the work of Huang {\em et al.} \cite{HLL}, Gasull {\em et al.} \cite{GYZ}, Llibre {\em et al.}  \cite{LYZ} or Carbonell and Llibre \cite{CL1989} . This question has also been approached using bifurcations by, for instance, Benterki and Llibre \cite{BL} or \cite{GYZ}.
Problems with symmetry appear in \'Alvarez {\em et al.} \cite{AGP}.
Our references do not pretend to be comprehensive. The reader can find further interesting work by looking at the references within those we provide.

We are, of course, also concerned in establishing an upper bound for the number of limit cycles. However, when no limit cycle exists, we take a different route and address the  question of the existence of policycles (sometimes called  heteroclinic cycles)  and the number of equilibria in them.

As many authors before us,  we are concerned with polynomial vector fields with homogeneous nonlinearities: vector fields of the form $\dot{X}=\lambda X + Q(X)$, where $X \in \RR^2$ and $Q$ is homogeneous.
However, our focus is on the special case where the nonlinear part is contracting, when Field's Invariant Sphere Theorem  \cite[Theorem 5.1, Theorem~\ref{IST} below]{Field89} guarantees the existence of an invariant circle for the dynamics.
Contracting nonlinearities occur quite naturally in some settings.
We provide a classification of the global dynamics for all such problems and obtain a complete invariant for the dynamics, including the behaviour  at infinity.

Cima and Llibre in  \cite{AnnaL1990Bond} define bounded vector fields in the plane and provide a classification of their behaviour at infinity.
Since vector fields with contracting nonlinearities are bounded in their sense, our results complement theirs by extending the classification globally.

The classification is then used to address some classes of examples.
We start with definite nonlinearities, that have been addressed by Gasull {\em et al} \cite{GLS1987}.
When the nonlinear part of the vector field is a contracting cubic, we are able to provide the full list of global phase portraits by making use of the results in Cima and Llibre \cite{AnnaL1990}.
If  the vector field is additionally $\Ztwo \oplus \Ztwo$-equivariant, we provide a complete description of the global planar dynamics, including the study of stability and bifurcation of equilibria.

\subsection*{Structure of the article}
In the next section we  establish some notation and state some results that will be used.
A normal form for planar contracting vector fields and sufficient conditions for a planar vector field to be contracting are obtained in Section~\ref {secD2}.
Dynamics  is discussed in Section~\ref{secGD2}  for the restriction to the invariant circle and globally in Section~\ref{subsec:classification}, where we  also obtain a complete invariant for the dynamics and from it a complete classification of this type of vector fields.
This is used in the remainder of the article to obtain a complete description of some families of examples: definite nonlinearities in Section~\ref{sec:exDefinite}; cubic nonlinearities in Section~\ref{sec:cubic};
$\ZZ_2\oplus\ZZ_2$-equivariant nonlinearities as special cases   in Subsections \ref{subsec:Z2dyn} and in \ref{sub:Z2}.

\section{Preliminary results and notation }\label{sec:prelim}
In this article we are concerned with the differential equation
\begin{equation} \label{eqR2contracting}
\left\{\begin{array}{lcl}
\dot{x}&=& \lambda x +Q_1(x,y)  \\
\mbox{} & \mbox{} & \\
\dot{y} &=& \lambda y +Q_2(x,y)
\end{array}  \right.
 \qquad\mbox{with}\qquad
\lambda >0
\end{equation}
where the $Q_i$, $i=1,2$  are homogeneous nonzero polynomials
%where $\lambda \neq 0$ and each $Q_i$, $i=1,2$  is a homogeneous polynomial
 of the same degree  $n>1$ and $\xy\in\plano$.
 % (in particular none of them is the zero polynomial)}.
 %$n$.
 %In this instance, w
 We define $Q=(Q_1, Q_2)$ and say it is a homogeneous polynomial of degree $n$.
  The origin of such a system is an  {\em unstable star node}, a node with equal and positive eigenvalues.
  
  For $\lambda<0$ the origin is an attracting star node and the dynamics corresponds to the equation with $Q$ replaced by $-Q$ and reversed time orientation.
  
  We recall some elementary notions in (equivariant) dynamical systems. The standard reference is the book \cite{G&S}.
We say that the dynamical system described by an ordinary differential equation $\dot{X} = f(X)$, $X \in \RR^n$ is {\em equivariant} under the action of a compact Lie group $\Gamma$ if
$$
f(\gamma .X) = \gamma .f(X)
$$
for all $X \in \RR^n$ and $\gamma \in \Gamma$. An {\em equilibrium} of $\dot{X} = f(X)$ is a solution of $f(X)=0$, the form \eqref{eqR2contracting} implies that at least the origin is an equilibrium. A {\em limit cycle} is an isolated periodic orbit. A {\em policycle} is the cyclic union of finitely many equilibria and  trajectories connecting them.

  Let $\langle , \rangle$ denote the  inner product and   $|| .  ||$ the norm in $\plano$, and let
 $P^{d}(\plano,\plano )$ be
%$\PP^d=P^{d}(V, V )$ 
the vector space of homogeneous polynomial maps of degree $d$ from $\plano$ to itself.
Denote by $P^{d+1}(\plano, \RR )$ 
%$\VV^{d+1}=P^{d+1}(V, \RR )$ 
the vector space of homogeneous polynomial maps of degree $d+1$ from $\plano$ to $\RR$ and let $X\in \plano$.
Consider the linear map:
\begin{equation}\label{eq:M}
\MM:P^{d}(\plano, \plano )\seta P^{d+1}(\plano, \RR )
\qquad
\MM Q(X)=\langle X,Q(X)\rangle .
\end{equation}
A polynomial $Q \in P^{d}(\plano, \plano)$, $d> 1$
%$Q \in  P^{2p+1}(V, V )$, $p \geq 1$,
is said to be {\em contracting} if 
$$
\MM Q(X) < 0, \quad  \mbox{for all }\quad || X  || =1.
$$
%\langle Q(\bw ), \bw \rangle < 0, \quad  \mbox{for all }|| \bw  || =1.
%$$
It follows
 that polynomials of even degree are never contracting. It is also useful to recall that, since the polynomial is homogeneous, stating that the inequality in the definition of contracting holds on the unit sphere is equivalent to saying that it holds for any nonzero vector.
We will also need the linear map $\LL: P^{2p+1}(\RR^2, \RR^2 )\seta  P^{2p+2}(\RR^2, \RR )$, given by
\begin{equation}\label{eq:L}
\LL Q(X)=\langle X^\perp,Q(X)\rangle
\qquad\mbox{for}\qquad
X=(x,y)
\qquad\mbox{and}\qquad
X^\perp=(-y,x) .
\end{equation}

 For ease of reference we state next a two-dimensional version
% \footnote{\isl   Optei por enunciar o teorema s\'o em dimens\~ao 2, fica mais simples e mais f\'acil de bater certo.
%O problema era que a forma do $\mathcal{P}_Q(\bw )=Q(\bw )-\langle Q(\bw ), \bw \rangle \bw $ usada pelo Field \'e  a proje\ca o do campo sobre o espa\c co tangente da esfera.
%O resultado d\'a um campo de grau $2p+3$, ou seja, grau de $Q$ mais 2, n\~ao  homog\^eneo.
%N\~ao d\'a o campo em coordenadas polares que est\'a descrito em \eqref{eq:polar_n}, como n\'os dizemos na sec. \ref{secGD2}.
%Acontece que em dimens\~ao 2 fica tudo mais simplificado.}
 of the Invariant Sphere Theorem \cite[Theorem 5.1]{Field89}, which  we will use extensively.

\begin{theorem}[The Invariant Sphere Theorem] \label{IST}
Let $p \geq 1$ and suppose that $Q \in  P^{2p+1}(\plano, \plano )$ is contracting. 
Then, for every $\lambda > 0$, there exists a unique topological circle
 %n-dimensional sphere 
 $S(\lambda) \subset  \plano \setminus \{0\}$ which is invariant by the flow of \eqref{eqR2contracting}.
%\begin{equation}\label{eqR2contracting}
%  \dot{x}=\lambda x + Q(x).
%\end{equation}
%
Further,
\begin{itemize}
  \item [(a)] $S(\lambda)$ is globally attracting in the sense that every trajectory $(x(t),y(t))$ of \eqref{eqR2contracting}
with nonzero initial condition is asymptotic to $S(\lambda)$ as $t\to +\infty$.
  \item [(b)]  $S(\lambda)$ is embedded as a topological submanifold of $\plano$ and the bounded
component of $\plano \setminus S(\lambda)$ contains the origin.
  \item [(c)] The flow of \eqref{eqR2contracting} restricted to $S(\lambda)$ is topologically equivalent to the flow
of the phase
%\footnote{\isl A equa\ca o era $\dot\theta=r^{2p}g(\theta)$. Acho que podemos tirar  o $r^{2p}$ porque multiplicar por uma fun\ca o positiva preserva a equival\^encia topol\'ogica. O mesmo vale para a din\^amica no infinito, a seguir.} 
equation  $\dot\theta=g(\theta)$ where $g(\theta)=\LL Q(\cos\theta,\sin\theta)$.
%vector field $\mathcal{P}_Q(\bw ) = Q(\bw )-\langle Q(\bw ), \bw \rangle \bw $.
\end{itemize}
\end{theorem}

%For planar vector fields, we have $V=\RR^2$ and the invariant sphere is 1-dimensional.
The odd degree of the nonlinear part $Q$ in the statement of Theorem~\ref{IST} implies that the vector field is 
$\ZZ_2$-equivariant, where $\ZZ_2$ is generated by $-Id$.

We will use  the representation of \eqref{eqR2contracting} in polar coordinates $(x,y)=(r\cos\theta,r\sin\theta)$, with $(r,\theta)\in \mathbb{R}^+\times \mathcal{S}^1$. This is given by
\begin{equation}\label{eq:polar_n}
\left\{ \begin{array}{l}
\dot{r} = \lambda r + f(\theta) r^{2p+1} \\
\dot{\theta} = g(\theta) r^{2p}
\end{array} \right.
\qquad\mbox{with}\qquad
\begin{array}{l}
f(\theta)=\MM Q(\cos\theta,\sin\theta)\\
%=\cos\theta Q_1(\cos\theta,\sin\theta) +\sin\theta Q_2(\cos\theta,\sin\theta)\\
{}\\
g(\theta) =\LL Q(\cos\theta,\sin\theta)
%= \cos\theta Q_2(\cos\theta,\sin\theta) -\sin\theta Q_1(\cos\theta,\sin\theta)
\end{array}
\end{equation}
where $\LL$ and $\MM$ are the linear maps defined in \eqref{eq:M} and \eqref{eq:L}.

Let $\contr^{2p+1}\subset  P^{2p+1}(\RR^2,\RR^2 )$ denote the set of contracting polynomial vector fields.
Our aim is to describe the global dynamics of \eqref{eqR2contracting} for $Q\in\contr^{2p+1}$, $p\ge 1$,  including the behaviour at infinity using the Poincar\'e disc, a  compactification  of $\RR^2$ (see Chapter 5 of Dumortier {\em et al.} \cite{DLA}).
The plane $\RR^2$ is identified to a compact disc, with its boundary corresponding to infinity.
The disc is also identified to a hemisphere in the unit sphere  $\mathcal{S}^2\subset \RR^3$, covered by six charts
$U_i, V_i$, $i=1,2,3$.
 In the coordinates $(u,v)$ on any of the charts $v=0$ corresponds to the equator  $\mathcal{S}^1$ of the sphere, the {\em  circle at infinity} of the Poincar\'e disc.
A point with coordinates $(u,v)$, $u\ne 0$ in $U_1$ corresponds to the point with coordinates $(\tilde{u},\tilde{v})=\left({1}/{u},{v}/{u}\right)$ in $U_2$ and to the point with coordinates $(\hat{u},\hat{v})=\left({u}/{v},{1}/{v}\right)$ in $U_3$.
The dynamics of \eqref{eqR2contracting} in the charts $U_1$ and $U_2$  are given, respectively, by
\begin{equation} \label{eqPXUboth}
\left\{\begin{array}{lcl}
\dot{u} &=& Q_2(1,u)-uQ_1(1,u)\\
{} &{} & \\
\dot{v} &=& -\lambda v^{2p+1} -vQ_1(1,u)
\end{array}  \right.
\quad\mbox{and}\qquad
\left\{\begin{array}{lcl}
\dot{u} &=& Q_1(u,1)-uQ_2(u,1) \\
{} & {} & \\
\dot{v} &=&  -\lambda v^{2p+1} -vQ_2(u,1)
\end{array}  \right. 
%\quad\mbox{and}\qquad
%\left\{\begin{array}{lcl}
%\dot{u} &=& Q_1(u,v) \\
%{} & {} & \\
%\dot{v} &=& Q_2(u,v)
%\end{array}  \right. .
\end{equation}
and the expression on the chart $U_3$ is just \eqref{eqR2contracting} computed at $(x,y)=(u,v)$.
The expressions of the Poincar\'e compactification in the three remaining charts $V_j$ are the same as in $U_j$.

The dynamics  at infinity of \eqref{eqR2contracting} is thus given by the restriction of each one of the expressions in \eqref{eqPXUboth}  to the flow-invariant line $(u,0)$, since the second equation is trivially satisfied for $v=0$. An equilibrium at infinity of \eqref{eqR2contracting}  is an equilibrium  $(u,0)\in \mathcal{S}^1$ of one of the two equations.
We refer to it as an {\em infinite equilibrium}, by opposition
to {\em finite equilibria} $(u,v)$, $v\ne 0$.

The  dynamics of the restriction of \eqref{eqPXUboth} to the circle at infinity $(u,0)$  does not depend  on the linear part of \eqref{eqR2contracting}. Hence, it is equivalent to the dynamics of the phase equation  $\dot\theta=g(\theta)$ equivalent to that in \eqref{eq:polar_n}.
%$$
%\left\{\begin{array}{lcl}
%\dot{x} &=& Q_1(x,y)  \\
%\mbox{} & \mbox{} & \\
%\dot{y} &=& Q_2(x,y)
%\end{array}  \right.
%$$

\section{Contracting polynomial vector fields in dimension 2}\label{secD2}

The results in this section describe the homogeneous polynomial planar vector fields and provide conditions for ensuring these are contracting. 
In this way we obtain a description of vector fields \eqref{eqR2contracting} to which Theorem~\ref{IST} applies.

\begin{proposition}\label{prop:homogeneous_polynomial}
Any homogeneous polynomial vector field
$Q\xy=\left(Q_1\xy,Q_2\xy\right)$
in $\plano$ of degree $2p+1$
may be written in the form
\begin{equation}\label{formNonHomogeneous}
Q\xy=p_1(x^2,y^2)\left( x,0\right)+p_2(x^2,y^2)\left( 0,y\right)+
p_3(x^2,y^2)\left( y,0\right)+p_4(x^2,y^2)\left( 0,x\right)
\end{equation}
where $p_j(u,v)$, $j=1,\ldots,4$ are homogeneous polynomials of degree $p$.
\end{proposition}

\begin{proof}
Each vector monomial occurring in $Q$ has the form $a x^ky^\ell e_j$ where $e_j$ is the $j$-th vector of the canonical basis and  $k+\ell=2p+1$, hence in each case one of $k,\ell$ is odd and the other is even. Then $xp_1(x^2,y^2)$ is the sum of the vector monomials in $Q_1$ with odd $k$, and $yp_3(x^2,y^2)$ is the sum of those with odd $\ell$. Similarly, $yp_2(x^2,y^2)$ is the sum of the vector monomials in $Q_2$ with odd $\ell$, and $xp_4(x^2,y^2)$ is the sum of those with odd $k$.
\end{proof}

We call $p_1(x^2,y^2)\left( x,0\right)+p_2(x^2,y^2)\left( 0,y\right)$ the {\em symmetric part} of $Q$ and
$p_3(x^2,y^2)\left( y,0\right)+p_4(x^2,y^2)\left( 0,x\right)$ the {\em asymmetric part} of $Q$. We write $Q_s(x,y)$ for the symmetric part of $Q$ and
note that it is $\Ztwo\oplus\Ztwo$-equivariant, where $\Ztwo\oplus\Ztwo$ is the group generated by the maps 
$\xy\mapsto(-x,y)$ and $\xy\mapsto(x,-y)$.

\begin{proposition}\label{prop:ContractSufficient}
A homogeneous polynomial vector field $Q$ of degree $2p+1$  in $\plano$ is contracting if
for the polynomials in \eqref{formNonHomogeneous}, we have  for all $(u,v)\ne(0,0)$ with  $u\ge 0$, $v\ge 0$,
that one of the $p_j(u,v)<0$, $j=1,2$ and
$$
2\max_{j=1,2}\left\{p_j(u,v)\right\}< - \left|p_3(u,v)+ p_4(u,v) \right|.
$$
\end{proposition}

 Note that if  $p_1(u,v)<0$, then the second condition implies $p_2(u,v)<0$ and vice-versa.
\begin{proof}
We have $\MM Q\xy=\xy\cdot A(x^2,y^2)\cdot \xy^T$ where $A$
 is the symmetric matrix
$$
A(u,v)=\left( \begin{array}{cc}
p_1(u,v)&\dfrac{p_3(u,v)+p_4(u,v)}{2}\\
&\\
\dfrac{p_3(u,v)+p_4(u,v)}{2}&p_2(u,v)
\end{array}\right) .
$$
The polynomial $Q$ is contracting if for each $(u,v)$ the quadratic form $\xy\cdot  A(u,v)\cdot \xy^T$ is negative definite.
This holds
if and only if both eigenvalues of $A(u,v)$ are negative.
By Gershgorin's Theorem \cite[Section 2.7.3]{Demmel} the eigenvalues of $A$ lie in the union of the closed intervals  with centre at $p_j(u,v)$, $j=1,2$ and radius $\left|p_3(u,v)+p_4(u,v) \right|/2$.
The inequality implies that both these intervals are contained in the negative half line.
\end{proof}

\begin{proposition}\label{prop:ContractGeneral}
A homogeneous polynomial vector field $Q$ of degree $2p+1$  in $\plano$ is contracting if  for the polynomials in \eqref{formNonHomogeneous}, for all $(u,v)\ne(0,0)$ with  $u\ge 0$, $v\ge 0$,
one of the $p_j(u,v)<0$, $j=1,2$ and
 \begin{equation}\label{Ca}
 4p_1(u,v)p_2(u,v)>(p_3(u,v)+p_4(u,v))^2.
 \end{equation}
\end{proposition}

\begin{proof}
The eigenvalues of the matrix  $A$
of the proof of Proposition~\ref{prop:ContractSufficient} are negative if and only if
$\tr A(u,v)<0$ and $\det A(u,v)>0$.
In case \eqref{Ca} holds then the hypothesis on the sign of one $p_j(u,v)$, $j=1,2$ implies that both $p_j(u,v)<0$, $j=1,2$ and hence that $\tr A=p_1(u,v)+p_2(u,v)<0$.
The result follows from
$\det A=p_1(u,v)p_2(u,v)-\left(p_3(u,v)+p_4(u,v)\right)^2/4$.
\end{proof}

The conditions in Propositions~\ref{prop:ContractSufficient} and \ref{prop:ContractGeneral} are not necessary.
A simple example is the symmetric vector field with $p_1\xy=y^2-x^2$, $p_2\xy=-2x^2-y^2$, $p_3\xy=p_4\xy=0$, 
for which $\MM Q\xy=-(x^4+y^4+x^2y^2)<0$ for $\xy\ne(0,0)$, but $p_1\xy=0$ for $x=\pm y$.

\begin{corollary}
 If $Q$ is a polynomial vector field satisfying the hypothesis of either Proposition~\ref{prop:ContractSufficient} or \ref{prop:ContractGeneral} then its symmetric part $Q_s$  is also contracting.
\end{corollary}

\begin{proof}
We have $Q_s\xy= p_1(x^2,y^2)\left( x,0\right)+p_2(x^2,y^2)\left( 0,y\right)$.
Hence it  follows  that both $p_1(x^2,y^2)$ and $p_2(x^2,y^2)$ are negative.
Since in this case $\MM Q\xy=\xy \cdot D\cdot \xy^T$ where
$$
D=\left( \begin{array}{cc}
p_1(x^2,y^2)&0\\0&p_2(x^2,y^2)
\end{array}\right)
$$
the definition of a contraction is satisfied.
\end{proof}

\section{Dynamics on the invariant circle}\label{secGD2}
%\section{Global dynamics for contracting nonlinearities.}\label{secGD2}

The hypothesis of contracting homogeneous nonlinearities in the vector field given by \eqref{eqR2contracting}, allows us to apply Theorem~\ref{IST}, guaranteeing the existence of a globally attracting invariant circle.
Observe that, from the expression in polar coordinates \eqref{eq:polar_n}, the homogeneous polynomial vector field $Q$ is contracting if and only if $f(\theta)<0$ for all $\theta$.

The form of the phase vector field on the invariant circle $\mathcal{S}^1\subset \mathbb{R}^2$ in Theorem~\ref{IST} is $\dot{\theta}=g(\theta)=\LL Q(\cos\theta,\sin\theta)$.
It determines the same dynamics as the expression  \eqref{eq:polar_n} for $\dot\theta$  in polar coordinates,
%$\bw =(w_1,w_2)=(\cos\theta,\sin\theta)$, 
since they differ by a positive function $r^{2p}$.
%and is given by
%$$
%\dot{\theta}=g(\theta)=\LL Q(\bw).
%$$
It follows that the dynamics on the invariant circle coincides with the dynamics on the circle at infinity.
We explore this in the following results, starting with three lemmas that are immediate.
These results are strongly related to \cite{ACL2021} and \cite{BLS2013}.

\begin{lemma}\label{lem:InvCycle}
Assume that $Q$ in \eqref{eqR2contracting} is contracting and homogeneous.
The invariant circle that exists for the dynamics of  \eqref{eqR2contracting} is an attracting limit cycle if and only if $g(\theta) \neq 0$ for all $\theta \in [0,2\pi)$.
Moreover, in this case the invariant circle is the curve
%\footnote{\isl Aqui estava ``the invariant circle is the graph of the map $r(\theta)=\sqrt[n-1]{-\lambda/f(\theta)}$'' mas estava errado porque esta \'e a equa\c c\~ao para $\dot r=0$, mas em geral $r$ varia nesta curva, a menos que $f$ seja constante.} 
$\LL Q\xy=0$
and another periodic orbit exists at infinity.
\end{lemma}

\begin{figure}[htbp]
\begin{center}
\parbox{25mm}{
\includegraphics[width=25mm]{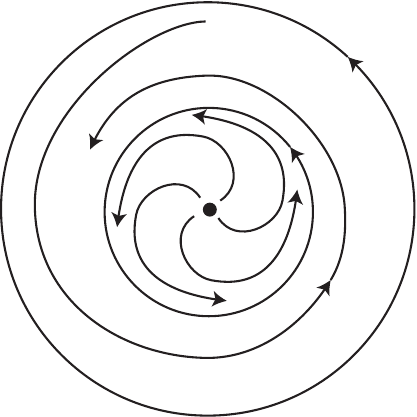}\\
(a) $g(\theta)\neq 0$\\
limit cycle
}
\qquad
\parbox{30mm}{
\includegraphics[width=25mm]{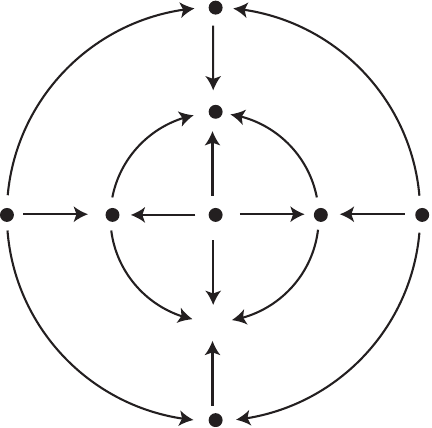}\\
(b) $g(\theta)= 0$ at\\
finitely many\\
 equilibria
}
\qquad
\parbox{30mm}{
\includegraphics[width=25mm]{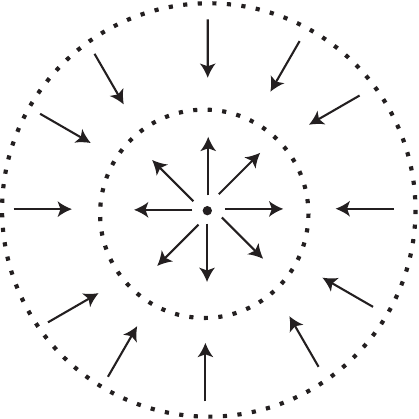}\\
(c) $g(\theta) \equiv 0$\\
  infinitely many equilibria
}
\end{center}
\caption{Global planar dynamics with star nodes as described in Proposition~\ref{prop:global_dynamics pi} and Lemmas~\ref{lem:InvCycle} -- \ref{lem:Cont}.}\label{fig:global_dynamics}
\end{figure}

\begin{proof}
There are no equilibria on invariant circle and on the circle at infinity since the phase equation has no zeros, hence both circles are limit cycles. The form of the invariant circle follows from the invariance of this curve established in  \cite[Theorem 1 (a)]{BLS2013}, see also Figure~\ref{fig:global_dynamics}~(a).
\end{proof}

\begin{lemma}\label{lem:policycle}
Assume that $Q$ in \eqref{eqR2contracting} is contracting and homogeneous.
The invariant circle that exists for the dynamics of  \eqref{eqR2contracting} is an attracting policycle if and only if $g(\theta) = 0$ for a finite number of $\theta \in [0,2\pi)$.
Moreover, in this case another policycle exists at infinity.
\end{lemma}

\begin{proof}
Both the  invariant circle and  the circle at infinity contain equilibria, hence they must be policycles as in Figure~\ref{fig:global_dynamics}~(b).
\end{proof}

\begin{lemma}\label{lem:Cont}
Assume that $Q$ in \eqref{eqR2contracting} is contracting and homogeneous of degree $2p+1$.
The invariant circle that exists for the dynamics of  \eqref{eqR2contracting} is a continuum of equilibria if and only if $g(\theta)= 0$ for all $\theta \in [0,2\pi)$.
Moreover, in this case the invariant circle is the graph of the map $r(\theta)=\sqrt[2p]{-\lambda/f(\theta)}$
and the circle at infinity is also a continuum of equilibria.
\end{lemma}

\begin{proof}
The phase equation  being identically zero, both the invariant circle and  the circle at infinity consist of equilibria. In polar coordinates, finite equilibria must also satisfy $\dot r=0$ and this provides the equation for the invariant circle. Phase portrait in Figure~\ref{fig:global_dynamics}~(c).
\end{proof}

\begin{proposition} \label{prop:global_dynamics pi}
Consider \eqref{eqR2contracting} with  $Q(x,y)\in\contr^{2p+1}$  a contracting polynomial vector field in the form given by  \eqref{formNonHomogeneous}.
Then:
\begin{enumerate}
\renewcommand{\theenumi}{(\alph{enumi})}
\renewcommand{\labelenumi}{{\theenumi}}
 \item\label{item:globala}
 If $p_3(u,v)p_4(u,v)<0$ and $-{4}p_3(u,v)p_4(u,v)>(p_2(u,v)-p_1(u,v))^2$ the invariant circle is a limit cycle.
  \item\label{item:globalb}
  If $p_3(0,1)p_4(1,0)\geq 0$ the invariant circle is either a policycle with at most $4(p+1)$ equilibria or a continuum of equilibria.
  Moreover, if  $p_3(0,1)p_4(1,0)> 0$ then the invariant circle is a policycle.
  \item\label{item:globalc}
  If $p_3(u,v)\equiv p_4(u,v)\equiv 0$ and  $p_1(u,v)\equiv p_2(u,v)$, the invariant circle is a  continuum of equilibria.
\end{enumerate}
\end{proposition}

Note that in case \ref{item:globalc} the equations are $\ZZ_2\oplus \ZZ_2$-equivariant, this property will be explored further in Subsections~\ref{subsec:Z2dyn} and \ref{sub:Z2}.
We illustrate in Figure~\ref{fig:global_dynamics} the possibilities described in Proposition~\ref{prop:global_dynamics pi}.

\begin{proof}
Using \eqref{formNonHomogeneous} we can write $\LL Q\xy=\xy\cdot B(x^2,y^2) \cdot\xy^T$ where
$$
B(x^2,y^2)=\begin{pmatrix}p_4(x^2,y^2) &\left(p_2(x^2,y^2)-p_1(x^2,y^2)\right)/2\\
\left(p_2(x^2,y^2)-p_1(x^2,y^2)\right)/2&-p_3(x^2,y^2)
\end{pmatrix}.
$$
If $p_3(u,v)p_4(u,v)<0$ and $-4p_3(u,v)p_4(u,v)>(p_2(u,v)-p_1(u,v))^2$  then $\det B>0$.
Hence if $p_4(u,v)>0$ then $g(\theta)=\LL Q(\cos\theta,\sin\theta)>0$ for all $\theta$,
with $g(\theta)<0$ provided $p_4(u,v)<0$. 
In both cases $g(\theta)\neq 0$ for all $\theta \in[0,2\pi]$ and item~\ref{item:globala} holds by Lemma~\ref{lem:InvCycle}.

If either $p_3(0,1)= 0$ or $p_4(1,0) = 0$ then trivially $g(\theta)$ vanishes on one of the axes and one of Lemmas~\ref{lem:policycle} and \ref{lem:Cont} holds. 
Assume $p_3(0,1)p_4(1,0) \neq 0$. 
Since $g(0)= p_4(1,0)$ and $g(\frac{\pi}{2}) = -p_3(0,1)$, then 
 $g(\theta)\not\equiv 0$ and thus the invariant circle is not a continuum of equilibria.
Moreover, in this case
$g(\theta)$ changes sign in the interval $(0,\pi/2)$. 
Therefore, since $g$ is continuous, there must be at least one $\theta^* \in (0,\pi/2)$ for which $g(\theta^*)=0$. 
Hence Lemma~\ref{lem:policycle}  applies and \eqref{eqR2contracting} has a policycle, establishing \ref{item:globalb}.

Item~\ref{item:globalc} is an immediate consequence of Lemma~\ref{lem:Cont}. 
\end{proof}

The next example illustrates a situation not accounted for by Proposition~\ref{prop:global_dynamics pi}.

\begin{Ex}
The family of vector fields studied by Boukoucha \cite{B2017} is such that $p_3(0,1)p_4(1,0)<0$ and a limit cycle exists. When $n=1$ in \cite{B2017}, we obtain
\begin{eqnarray*}
p_1(x^2,y^2) = -\beta a x^2 -(\beta a + \alpha b) y^2 & \mbox{and } & p_3(x^2,y^2) = (\alpha a+\beta b)x^2+\alpha a y^2 \\
p_2(x^2,y^2) = (\alpha b-\beta a)x^2 - \beta a y^2 & \mbox{and } & p_4(x^2,y^2) = (\beta b - \alpha a) y^2-\alpha a x^2
\end{eqnarray*}
for real constants $\alpha, \beta, a, b$.
Then $p_3(0,1)p_4(1,0) = -\alpha^2 a^2 <0$.
If $b=0$,  then $p_1\equiv p_2$ and $p_3+p_4\equiv 0$.
Hence if $\beta a>0$ then  $Q$ is contracting  by Proposition~\ref{prop:ContractSufficient}.
Moreover, $p_3p_4(u,v)=-\alpha^2 a^2(u+v)^2<0$ and this example satisfies the conditions in Proposition~\ref{prop:global_dynamics pi} (a).

Another choice of parameters for which $Q$ is contracting  is  $\alpha=0$, $\beta a>0$ and $4a^2-b^2>0$, this time  by Proposition~\ref{prop:ContractGeneral}.
In this case $p_3\xy=b\beta x^2$, $p_4\xy=b\beta y^2$ and $p_1\xy=p_2\xy=-\beta a(x^2+y^2)$.
Therefore with this choice of parameters and if $b\beta\ne 0$ the example does not satisfy the conditions in Proposition~\ref{prop:global_dynamics pi} (a) and yet the invariant circle is a limit cycle.
\end{Ex}

\begin{corollary}\label{cor:stability}
If  $Q(x,y)$ is  a contracting polynomial vector field for which \eqref{eqR2contracting} has a finite number of equilibria then:
\begin{enumerate}
\renewcommand{\theenumi}{(\alph{enumi})}
\renewcommand{\labelenumi}{{\theenumi}}
\item\label{item:hyperbolic}
if all the equilibria of \eqref{eqR2contracting} are hyperbolic, then the number of  equilibria away from the origin is a multiple of 4 and they alternate as sinks and saddles;
\item\label{item:non-hyperbolic}
all the equilibria of \eqref{eqR2contracting} away from the origin are either sinks or saddles (possibly non-hyperbolic) or saddle-nodes;
\item\label{item:alternate}
the equilibria that are sinks and saddles appear at alternating positions in the policycle.
\end{enumerate} 
\end{corollary}

\begin{Ex} \label{ex:n=5} The following vector field illustrates the global dynamics given in Proposition~\ref{prop:global_dynamics pi}~\ref{item:globalb} when the nonlinearity is of degree $2p+1=5$ (see Figure \ref{fig: prop stability_contracting}):
$$
\left\{\begin{array}{lcl}
\dot{x} &=& \lambda x -x(x^4+x^2y^2+y^4)-y(-x^4+x^2y^2)\\
\dot{y} &=& \lambda y -y(x^4+x^2y^2+y^4)-x(x^2y^2-y^4)
\end{array}  \right.
$$

\begin{figure}[hht]
 \begin{center}
  \includegraphics[width=40mm]{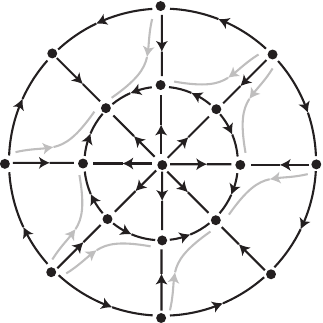}
  \end{center}
  \caption{Phase portrait of Example~\ref{ex:n=5}}\label{fig: prop stability_contracting}
\end{figure}

It follows by Proposition~\ref{prop:ContractGeneral} that the nonlinear part $Q\xy$ of this example is contracting because, for all $(x,y)\neq (0,0)$
%$Q(x,y)=(-x(x^4+x^2y^2+y^4)-y(-x^4+x^2y^2), -y(x^4+x^2y^2+y^4)-x(x^2y^2-y^4))$ is contracting provided that  
$$
p_1(x^2,y^2)=p_2(x^2,y^2)=-(x^4+x^2y^2+y^2)<0
$$ 
and 
$$
4p_1(x^2,y^2)^2-(p_3(x^2,y^2)+p_4(x^2,y^2))^2=3(x^2-y^2)^2+9x^4y^4+6x^2y^2(x^2-y^2)^2>0 .
$$
Then 
$$
g(\theta)=2\cos^2{\theta}\sin^2{\theta}\cos{(2\theta)}=\frac{1}{2}\sin^2{(2\theta)}\cos{(2\theta)}
$$
and
$$
g'(\theta)=\sin{(2\theta)}[2\cos^2{(2\theta)}-\sin^2{(2\theta)}].
$$
Hence, the four infinite equilibria on the axes are of saddle-node type and $\theta=\frac{\pi}{4}, \frac{3\pi}{4}, \frac{5\pi}{4}, \frac{7\pi}{4}$ are repellor, saddle, repellor and saddle, respectively. 
%attractor, saddle, attractor and saddle, respectively. 
Moreover, $p_3(0,1)=p_4(1,0)=0$ and the system has a policycle with the total number
%\footnote{\isl Aqui estava ``total number of equilibria away from the origin equal to $12$, the maximum predicted by Proposition~\ref{prop:global_dynamics pi}\ref{item:globalb}.'' Parece que n\~ao sabemos contar!}
 of equilibria away from the origin equal to 8.
\end{Ex}

\subsection{Special case:  $\ZZ_2\oplus\ZZ_2$ equivariant nonlinearity}\label{subsec:Z2dyn}
If the vector field \eqref{eqR2contracting} has $\Ztwo\oplus\Ztwo$ symmetry then
$Q$ has the form $Q\xy= p_1(x^2,y^2)\left( x,0\right)+p_2(x^2,y^2)\left( 0,y\right)$ and
 we may say more about the dynamics on the invariant circle.
In this case if $Q$ has degree $2p+1$ we may write 
$$
p_j\xy=\sum_k a_{jk}(x^2)^{p-k}(y^2)^{k}
\qquad j=1,2 .
$$

\begin{lemma}[Infinitely many equilibria]\label{corContinuum}
Let $Q$ be a $\Ztwo\oplus\Ztwo$-equivariant  contracting homogeneous polynomial  vector field  and suppose  $\lambda>0$. Then the invariant circle of \eqref{eqR2contracting} consists entirely of non hyperbolic equilibria if and only if $p_1(x^2,y^2)\equiv p_2(x^2,y^2)$.
\end{lemma}
\begin{proof}
If $p_1(x^2,y^2)\equiv p_2(x^2,y^2)$ then all points in the curve $\lambda = -p_1(x^2,y^2)$ are equilibria.
Conversely, all the points in the invariant circle are equilibria if and only if 
$\LL Q\xy=xy\left(p_2(x^2,y^2)-p_1(x^2,y^2)\right)\equiv 0$.
The equilibria are not hyperbolic since they form a continuum.
\end{proof}

When there are finitely many equilibria we use the polar form \eqref{eq:polar_n} and write
$$
g(\theta)=\LL Q(\cos\theta,\sin\theta)=\frac{1}{2} \sin(2\theta)\left[p_2\left(\frac{1+\zeta}{2},\frac{1-\zeta}{2}\right) -p_1\left(\frac{1+\zeta}{2},\frac{1-\zeta}{2}\right)\right] =\frac{1}{2} \sin(2\theta)q(\zeta)
$$
where $\zeta=\cos(2\theta)$.
Denote $\left(\xi_1,\xi_2\right)=\left(\dfrac{1+\zeta}{2},\dfrac{1-\zeta}{2}\right)$.
Then
$$
\begin{array}{lcl}
g'(\theta)&=&\cos(2\theta)\left[p_2\left(\xi_1,\xi_2\right)-p_1\left(\xi_1,\xi_2\right) \right]+\\
&&\dfrac{1}{2} \sin(2\theta)
\left[\dfrac{dp_2}{d\xi_2}\left(\xi_1,\xi_2\right)-\dfrac{dp_1}{d\xi_2}\left(\xi_1,\xi_2\right)+
\dfrac{dp_1}{d\xi_1}\left(\xi_1,\xi_2\right)-\dfrac{dp_2}{d\xi_1}\left(\xi_1,\xi_2\right)\right] .
\end{array}
$$

\begin{corollary}[Equilibria on the axes]\label{corEqAxes}
If $a_{10}-a_{20}\ne 0$ then $\theta=0$ and $\theta=\pi$ are hyperbolic equilibria; otherwise they are non hyperbolic.\\
If $a_{1p}-a_{2p}\ne 0$ then $\theta=\pi/2$ and $\theta=3\pi/2$ are hyperbolic equilibria; otherwise they are non hyperbolic.
\end{corollary}
\begin{proof}
On the horizontal axis $g'(\theta)=p_2(1,0)-p_1(1,0)=a_{10}-a_{20}$.
On the vertical axis $g'(\theta)=p_1(0,1)-p_2(0,1)=-a_{1p}+a_{2p}$.
\end{proof}

\begin{corollary}[Equilibria outside the axes]\label{corEqOutAxes}
Equilibria with $\theta\ne n\pi/2$, $n\in\ZZ$, are hyperbolic if and only if $\dfrac{dp_1}{d\theta}\ne\dfrac{dp_2}{d\theta}$.
\end{corollary}
\begin{proof}
In this case both $\cos(2\theta)\ne 0$, $\sin(2\theta)\ne 0$ and $p_1=p_2$.
Hence $g'(\theta)\ne 0$ if and only if
$$
\dfrac{dp_1}{d\theta}-\dfrac{dp_2}{d\theta}=\frac{dp_1}{d\xi_2}-\frac{dp_2}{d\xi_2}+\frac{dp_2}{d\xi_1}-\frac{dp_1}{d\xi_1}\ne 0.
$$
\end{proof}

\section{Global dynamics and classification}\label{subsec:classification}
Next we focus on the different possibilities for the dynamics of  \eqref{eqR2contracting} when the nonlinear part is a contracting homogeneous polynomial.
We classify the possible dynamical behaviour, up to a global planar homeomorphism that maps trajectories to trajectories, preserving the time orientation in each trajectory, plus a global rescaling of time.
This induces an equivalence relation on the set $\contr^{2p+1}$ of contracting homogeneous polynomial vector fields in $\plano$ of degree $2p+1$.
Given $Q_a,Q_b\in\contr^{2p+1}$ we indicate this equivalence relation as  $Q_a\sim Q_b$.

Since the set of positive definite polynomials is an open half cone in $P^{2p+2}(\plano,\RR)$ then its inverse image $\contr^{2p+1}\subset P^{2p+1}(\plano,\plano)$ under the linear map $\MM$ defined in \eqref{eq:M} is also  an open half cone in $P^{2p+1}(\plano,\plano)$.
The next result shows that  $\LL\left(\contr^{2p+1} \right)= P^{2p+2}(\plano,\RR)$ where $\LL$ is the  linear map defined in \eqref{eq:L} that generates the phase vector field.

\begin{theorem}\label{teo:contract}
Given a homogeneous polynomial $q\xy\in P^{2p+2}(\plano,\RR)$ of degree $2(p+1)$ there is a contracting homogeneous polynomial vector field $Q\xy\in\contr^{2p+1}$ for which $\LL Q\xy=q\xy$.
\end{theorem}
\begin{proof}
Write $q\xy=x^2b_1(x^2,y^2)+xyb_2(x^2,y^2)+y^2b_3(x^2,y^2)$ where $b_j(u,v)$ are homogeneous of degree $p$.
Let $Q$ be the vector field of
the form \eqref{formNonHomogeneous} in Proposition~\ref{prop:homogeneous_polynomial}  where,
for some $K>0$ to be determined, the $p_j$ are
$$
p_1(u,v)= -K(u^p+v^p)\quad
p_2(u,v)=b_2(u,v)+p_1(u,v)\quad
p_3(u,v)=-b_3(u,v)\quad
p_4(u,v)=b_1(u,v).
$$
Then $\LL Q\xy=x^2p_4(x^2,y^2)+xy\left[p_2(x^2,y^2)-p_1(x^2,y^2)\right]-y^2p_3(x^2,y^2)=q\xy$.

We want to choose $K$ so that the $p_j$ satisfy the conditions of Proposition~\ref{prop:ContractGeneral}.
Since $K>0$ then 
$$
\max_{t\in [0,\pi/2]}p_1(\cos t,\sin t)=-2^{1-p/2}K<0,
$$
hence $p_1(u,v)<0$ for $u\ge 0$, $v\ge 0$, $(u,v)\ne(0,0)$.
It remains to 
%ensure that \eqref{Ca} holds for all $(u,v)$ with $u\ge 0$, $v\ge 0$.
%With the choices above  we intend to 
find $K>0$ such that \eqref{Ca} holds for all $(u,v)$ with $u\ge 0$, $v\ge 0$, i.e., such that for $(u,v)=(x^2,y^2)$ we have:
$$
4K^2(u^p+v^p)^2-4K(u^p+v^p)b_2(u,v)- \left[b_1(u,v)-b_3(u,v) \right]^2>0
% 4 p_1^2(u,v)+4b_2(u,v) p_1(u,v)- \left[b_1(u,v)-b_3(u,v) \right]^2>0
\quad
\forall u\ge 0, v\ge 0.
$$
Since $D(u,v)=b_2^2(u,v)+(b_1(u,v)-b_3(u,v))^2\ge 0$, then if we find $K$ such that $2p_1\le-b_2-\sqrt{D}$  for  $u\ge 0$, $v\ge 0$, $(u,v)\ne(0,0)$ it will follow that $Q$ is contracting.
Let $M$ satisfy $M\le (-b_2-\sqrt{D})/2$ for $(u,v)=(\cos^2 t,\sin^2 t)$ with $0\le t\le \pi/2$.
Since $p_1$ and the $b_j$ are homogeneous  of the same degree then by taking $-K\le M/2^{1-p/2}$ the result is proved.
\end{proof}

We establish in this section that the global dynamics of  \eqref{eqR2contracting} for $Q\in \contr^{2p+1}$ is completely determined by $g(\theta)=\LL Q(\cos\theta,\sin\theta)$.
This feature  allows us to have a complete classification of vector fields in $\contr^{2p+1}$ from the point of view of the dynamics of \eqref{eqR2contracting},
by describing the equivalence relation induced by $\sim$ in the set $ P^{2p+2}(\plano,\RR)$ 

The natural classification in $ P^{2p+2}(\plano,\RR)$ is to allow linear changes of coordinates and multiplication by a nonzero constant,
that we will take to be always positive in order to preserve stability, as  discussed below.
This classification has good properties with respect to the topology induced in $ P^{2p+2}(\plano,\RR)$ by identifying the coefficients in the polynomials to points in $\RR^{2p+3}$.
In particular, it creates a Whitney stratification of $ P^{2p+2}(\plano,\RR)$.
It also translates well to $\contr^{2p+1}$ respecting the  dynamics in the invariant circle,  as the next simple result shows.

\begin{lemma}\label{lem:linearChange}
If $L:\plano\seta\plano$ is an invertible linear map and $\GG(X)=\LL Q(X)$ then the change of coordinates $L\widetilde{X}=X$ transforms  \eqref{eqR2contracting} into an equation with 
$\LL \widetilde{Q}\left(\widetilde{X}\right)=\left({\det L}\right)\GG\left(L\widetilde{X}\right)$.
%$\LL \widetilde{Q}\left(\widetilde{X}\right)=\dfrac{1}{\det L}\GG\left(L\widetilde{X}\right)$.
\end{lemma}

\begin{proof}
The linear part $\lambda X$ of equation  \eqref{eqR2contracting} commutes with every linear map of $\RR^2$.
Therefore, the  change of coordinates transforms $\dot X=\lambda X+Q(X)$ into  
$\dot {\widetilde{X}}=\lambda \widetilde{X}+L^{-1}Q\left(L\widetilde{X}\right)$.
Writing $X^\perp=\left(PX^T\right)^T$ where $P=\begin{pmatrix}0&-1\\1&0\end{pmatrix}$ we get
$$
  \begin{array}{ll}
 \LL\left(L^{-1}QL \right)\left(\widetilde{X}\right)
&
 =\left\langle P\widetilde{X},L^{-1}Q
 \left(L\widetilde{X}\right) \right\rangle
 =\left\langle \left(L^{-1}\right)^T P\widetilde{X},
Q \left(L\widetilde{X}\right) \right\rangle\\
 &
 =\left({\det L}\right)\left\langle  PL\widetilde{X},
 Q\left(L\widetilde{X}\right) \right\rangle
 =\left({\det L}\right)\GG\left(L\widetilde{X}\right) 
%  =\dfrac{1}{\det L}\left\langle  PL\widetilde{X},
% Q\left(L\widetilde{X}\right) \right\rangle
% =\dfrac{1}{\det L}\GG\left(L\widetilde{X}\right) 
 \end{array}
 $$
 since by Cramer's rule $PL=\dfrac{1}{\det L} \left(L^{-1}\right)^T P$.
\end{proof}

Under the equivalence induced by $\sim$, the classification in $ P^{2p+2}(\plano,\RR)$ under linear changes of coordinates gives rise to moduli: parametrised families of polynomials that share the same geometry.
For instance in Cima \& Llibre's \cite{AnnaL1990}  classification of $\PP^1$, that we use in Section~\ref{sec:cubic} below,
the families \ref{Contract1},  \ref{Contract2} and  \ref{Contract3} all contain a parameter $\mu$ that does not have a qualitative meaning for the dynamics.
The moduli arise from the  position of the
 roots of the polynomial $\LL Q$ in the projective space $\projective$,
 since a linear map on the plane is determined by its value at two points, so a linear change of coordinates only controls the position of two roots.
Therefore, $\sim$ induces a coarser equivalence relation in $ P^{2p+2}(\plano,\RR)$, since a homeomorphism would not have this restriction.
This is addressed in the next definition.

\begin{figure}
\begin{center}
\includegraphics[width=80mm]{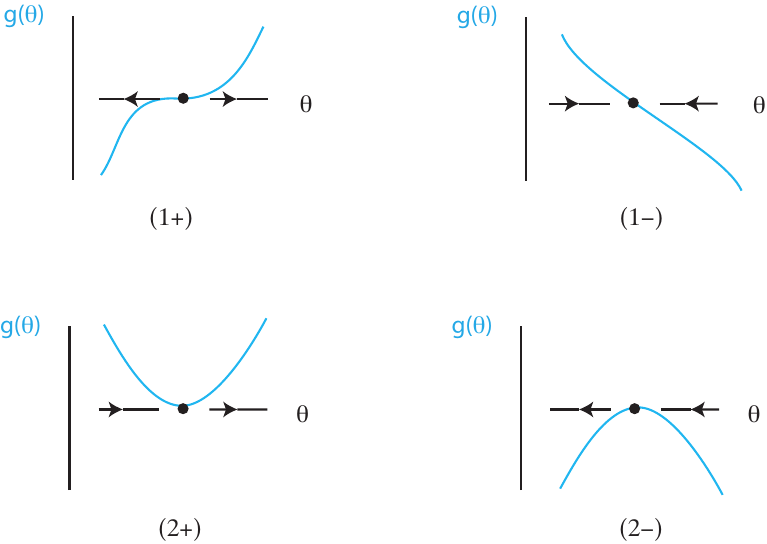}
\end{center}
\caption{Dynamics around the equilibria in the invariant circle and codes for the cyclic sequence. Codes (1+) and (1-) refer to repelling and attracting equilibria, respectively, codes (2+) and (2-) correspond the two possible orientations around a saddle-node.}\label{fig:localEquil}
\end{figure}

\begin{definition}\label{def:symbol}
The {\em symbol sequence} $\sigma(\GG)$associated to $\GG\in  P^{2p+2}(\plano,\RR)$ is a cyclic oriented list of the form $\sigma(\GG)=(j_1 s_1), (j_2 s_2),\ldots ,(j_\ell s_\ell)$ where $j_i\in \{1,2\}$ and $s_i=\pm$ obtained from the ordered set of zeros $0\le\theta_1<\theta_2<\cdots<\theta_\ell<\pi$ of $g(\theta)=\GG(\cos\theta,\sin\theta)$ as follows (see also Figure~\ref{fig:localEquil}):
\begin{enumerate}
%\item $j_i\in \{1,2\}$ is the modulo 2 multiplicity of $\theta_i$;
\item $j_i=1$ if the multiplicity of $\theta_i$ is odd and $j_i=2$ if it is even;
\item if $j_i=1$ then $s_i=+$  if $g(\theta)$ is increasing around $\theta_i$ and $s_i=-$  if $g$ is decreasing;
\item if $j_i=2$ then $s_i=+$  if $\theta_i$ is a local minimum of $g(\theta)$ and $s_i=-$  if $\theta_i$ is a local maximum;
\item if $g(\theta)\ne 0$ for all $\theta\in[0,\pi)$ then $\sigma(\GG)=\varnothing$;
\item if $\GG\xy\equiv 0$ then $\sigma(\GG)=\infty$.
\end{enumerate}
For $\sigma=(j_1 s_1), (j_2 s_2),\ldots,(j_\ell s_\ell)$, the {\em backward sequence}  is
$\bar\sigma=(j_\ell \tilde{s}_\ell), (j_{\ell-1} \tilde{s}_{\ell-1}),\ldots, (j_2 \tilde{s}_2),(j_1\tilde{s}_1)$, where $\tilde{s}_i=-s_i$ if $j_i=1$ and $\tilde{s}_i=s_i$ if $j_i=2$.
We identify $\sigma$ and $\bar\sigma$, and indicate this by
$\sigma\equiv \bar\sigma$.
\end{definition}

For instance, the symbol sequence for $\GG_1\xy=x^3y^2(x-y)$ is $\sigma(\GG_1)=(2+)(1-)(1+)$ corresponding to $\theta_1=0$ (double),
$\theta_2=\pi/4$ (simple) $\theta_3=\pi/2$ (triple).
For $\GG_2\xy=-x^2y^3(-x+y)=-\GG_1(y,x)$ the symbol sequence is $\sigma(\GG_2)=(1-)(1+)(2+)$ corresponding to $\theta_1=0$ (triple),
$\theta_2=\pi/4$ (simple) $\theta_3=\pi/2$ (double). 
Since the sequences are cyclic, they coincide.
Moreover, in this example $ \sigma(\GG_2)=\overline{\sigma(\GG_1)}$.

The sequence $\bar\sigma$ does not always coincide with $\sigma$.
An example is 
$$
\sigma=(2+),(1-),(2-),(1+),(1-),(1+)
\quad\mbox{with}\quad
\bar\sigma=(1-),(1+),(1-),(2-),(1+),(2+)
%\bar\sigma=(1-),(1+),(1-),(2-),(1+),(1-)
$$
where $\sigma=\sigma(\GG_1)$ for
$$
\GG_1\xy=x(a_1x-y)(a_2x-y)^2(a_3x-y)(a_4x-y)y^2
%\GG_1\xy=x(a_1x-y)(a_2x-y)^2(a_3x-y)(a_4x-y)(a_5x-y)y^2
$$
with $a_i=\tan(i\pi/5)$
%$a_i=\tan(i\pi/12)$
 and $\bar\sigma=\sigma(\GG_2)$ for 
 $\GG_2\xy=\GG_1(-x,y)$ (with the $(2+)$ moved to the beginning).
% $\GG_2\xy=-\GG_1(x,-y)$.

\begin{lemma}\label{lem:allowedSequences}
The symbol sequence $\sigma(\GG)=(j_1 s_1), (j_2 s_2),\ldots,(j_\ell s_\ell)$  of $\GG\in  P^{2p+2}(\plano,\RR)$ always satisfies the following restrictions:
\begin{enumerate}
\renewcommand{\theenumi}{(\alph{enumi})}
\renewcommand{\labelenumi}{{\theenumi}}
\item \label{item:0mod2}
$\sum_{i=1}^\ell j_i=0\pmod{2}$ ;
\item\label{item:alternating}
$(1+)$ and $(1-)$  occur in alternating sequences of sign + and -, even when the sequence is interrupted by one or more symbols $(2\pm)$;
\item\label{item:table}
If $s_i=+$ then $(j_{i+1},s_{i+1})\in\{(1-),(2+)\}$.
If $s_i=-$ then $(j_{i+1},s_{i+1})\in\{(1+),(2-)\}$.
\end{enumerate}
Moreover, if $\sigma$ satisfies these restrictions 
%of Lemma~\ref{lem:allowedSequences}, 
then $\bar\sigma$ also satisfies the same restrictions
\end{lemma}

\begin{proof}
Since the degree of  $\GG$ is even, restriction \ref{item:0mod2} follows.
The other two restrictions can be seen immediately from Figure~\ref{fig:localEquil}.
\end{proof}

The restriction \ref{item:alternating} corresponds to assertion \ref{item:alternate} in Corollary~\ref{cor:stability}.
Heteroclinic cycles occur for those $Q$ such that $\sigma(\LL Q)$ only contains one of the symbols $(2\pm)$.

\begin{proposition}\label{prop:LinearInvariant}
The symbol sequence $\sigma(\GG)$, under the identification $\equiv$, is  invariant  under linear changes of coordinates in $ P^{2p+2}(\plano,\RR)$.
\end{proposition}
\begin{proof}
Suppose $L:\plano\seta\plano$ is an invertible linear map and let 
$\GG_1\xy=\left({\det L}\right)\GG_2\circ L\xy$ 
%$\GG_1\xy=\dfrac{1}{\det L}\GG_2\circ L\xy$ 
with $g_j(\theta)=\GG_j(\cos\theta,\sin\theta)$, $j=1,2$.
Then $L$ maps the roots of $\GG_1$ in $\mathcal{S}^1$ into the  roots of $\GG_2$ with the same multiplicity.
Also there is a bijection $\varphi:\mathcal{S}^1\seta\mathcal{S}^1$ such that 
$g_2(\varphi(\theta))=\left({\det L}\right)g_1(\theta)$.
%$g_2(\varphi(\theta))=\dfrac{1}{\det L}g_1(\theta)$.
If $L$ preserves orientation, i.e. $\det L>0$, then the roots of $\GG_1$ and $\GG_2$ occur in the same order in $\mathcal{S}^1$.
The map $\varphi$ is monotonically increasing,  hence $\sigma(\GG_2)=\sigma(\GG_1)$.

If $L$ reverses orientation, i.e. $\det L<0$, then the roots of $\GG_1$ and $\GG_2$ occur in the opposite order in $\mathcal{S}^1$.
In this case the function $\varphi(\theta)$ is monotonically decreasing.
Hence,  if $g_1$ is a monotonically increasing (respectively, decreasing) function of $\theta\in[\theta_a,\theta_b]$ then $g_2$ is also a monotonically increasing (respectively, decreasing) function of $\tilde\theta=\varphi(\theta)$ for $\tilde\theta\in[\varphi(\theta_b),\varphi(\theta_a)]$.
Therefore $\sigma(\GG_2)=\overline{\sigma(\GG_1)}\equiv \sigma(\GG_1)$.
\end{proof}

In order to deal with the full equivalence relation $\sim$ in $\contr^{2p+1}$ we use
results of Neumann and O'Brien \cite{NewmannOBrien} for which we need to establish some terminology.
Let ${\mathbb D}$ be the Poincar\'e disc and let $\phi$ be the flow of  \eqref{eqR2contracting}.
Identifying each trajectory of  \eqref{eqR2contracting} to a point we obtain the {\em cell complex} $K(\phi)={\mathbb D}/\phi$, with  projection $\pi: {\mathbb D}\seta K(\phi)$ and some additional structure, as follows:
\begin{enumerate}
\renewcommand{\theenumi}{(\alph{enumi})}
\renewcommand{\labelenumi}{{\theenumi}}
\item\label{item:1d}
cells of dimension 1 correspond to {\em canonical regions}: open sets $A\subset{\mathbb D}$, homeomorphic to $\plano$ where the flow is equivalent to $\dot x=1$, $\dot y=0$;
\item\label{item:fiber}
cells $c$ of dimension 0 correspond to equilibria and separatrices of the flow and are initially classified by the dimension of the fibre $\pi^{-1}(c)$;
\item\label{item:order}
a partial order $<$ is defined on $K(\phi)$ as follows: separatrices in the boundary of canonical regions have the order induced by the flow; if $p$ is an equilibrium and $q$ is a point in a separatrix then 
if $p\in\alpha(q)$ then $\pi(p)<\pi(q)$,
if $p\in\omega(q)$ then $\pi(q)<\pi(p)$,
otherwise  $\pi(q)$ and $\pi(p)$ are not related.
\end{enumerate}

\begin{figure}
\parbox{50mm}{
\begin{center}
\includegraphics[width=35mm]{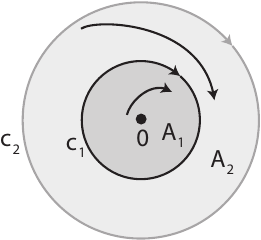}\\
(a)
\end{center}
}
\qquad\qquad
\parbox{60mm}{
\begin{center}
\includegraphics[width=19mm]{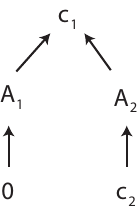}\\
(b)
\end{center}
}

\caption{(a) Phase portrait of \eqref{eqR2contracting}  when $g(\theta)\ne 0$. (b) Lattice for the partial order $<$ on $K(\phi)$ for this case. Arrows $x\to y$ indicate $\pi(x)<\pi(y)$. Cells $\pi(A_1)$ and $\pi(A_2)$ are 1-dimensional, for $A_1$ (respectively $A_2$) they are the projection of trajectories starting at the origin (respectively $c_2$) and ending at $c_1$. All the other cells are points. }\label{fig:limitCycle}
\end{figure}

\begin{figure}
\parbox{50mm}{
\begin{center}
\includegraphics[width=49mm]{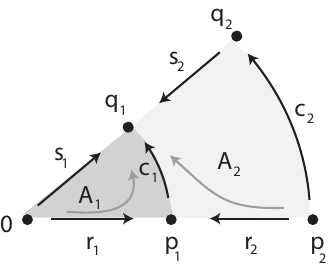}\\
(a)
\end{center}
}
\qquad\qquad
\parbox{60mm}{
\begin{center}
\includegraphics[width=59mm]{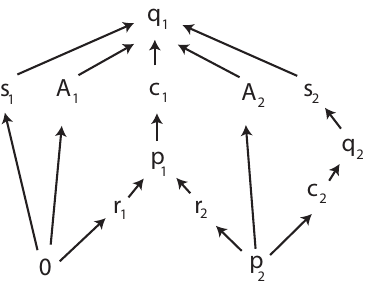}\\
(b)
\end{center}
}

\caption{(a) Phase portrait of \eqref{eqR2contracting}  on a sector. (b) Lattice for the partial order $<$ on $K(\phi)$ for the same sector: arrows $x\to y$ indicate $\pi(x)<\pi(y)$. Cells $\pi(A_1)$ and $\pi(A_2)$ are 1-dimensional, for $A_1$ (respectively $A_2$) they are the projection of trajectories starting at the origin (respectively $p_2$) and ending at $q_1$. All the other cells are points. }\label{fig:lattice}
\end{figure}

Examples are shown in Figures~\ref{fig:limitCycle} and \ref{fig:lattice}.

\begin{theorem}\label{prop:invariant}
The symbol sequence $\sigma(\LL Q)$, under the identification $\equiv$, is a complete invariant for the equivalence relation  $\sim$ in $\contr^{2p+1}$.
\end{theorem}

\begin{proof}
Let $\GG=\LL (Q)$ and $g(\theta)=\GG(\cos\theta,\sin\theta)$.
First suppose $g(\theta)\equiv 0$ or equivalently $\sigma(\GG)=\infty$. 
In this case, as in Lemma~\ref{lem:Cont}, all points in the invariant circle and in the circle at infinity are equilibria.
Apart from the origin all other trajectories are contained in rays, as in Figure~\ref{fig:global_dynamics} (c), hence  all $Q$ for which $\LL (Q)$ has this symbol sequence are equivalent.

The other simple case is $g(\theta)\ne 0$ for all $\theta$, as in Lemma~\ref{lem:InvCycle},  or equivalently $\sigma(\GG)=\varnothing$. 
The invariant circle and the circle at infinity are closed trajectories and the invariant circle attracts all finite trajectories not starting at the origin, by Theorem~\ref{IST}.
Apart from the origin all other trajectories are spirals, as in Figure~\ref{fig:global_dynamics} (a).
The cell complex consists of two  1-dimensional cells, two separatrices (the closed trajectories) giving rise to 0-dimensional cells with 1-dimensional fibre, and the equilibrium at the orgin yielding a  0-dimensional cell with 0-dimensional fibre, with the order shown in Figure~\ref{fig:limitCycle}.

Suppose now $g(\theta)=0$ at finitely many (and not zero) points, as in Lemma~\ref{lem:policycle}.
If $g(\theta_1)=0$ then, from the equation  \eqref{eq:polar_n}  in polar coordinates, it follows that the ray given by $\left\{(r,\theta_1)\ : \ r\ge 0 \right\} $ is flow-invariant.
Therefore two consecutive zeros  $\theta_1<\theta_2$ of $g$ define a flow-invariant {\em sector} 
  $$
  \left\{(r,\theta)\ : \ r\ge 0,\ \theta_1\le\theta\le \theta_2\right\} .
  $$
If $\theta_1<\theta_2<\theta_3$ are consecutive zeros of $g(\theta)$ we say the sector determined by $\theta_2$ and $\theta_3$ comes {\em after} the sector determined by $\theta_1$ and $\theta_2$.
The dynamics of  \eqref{eqR2contracting}  in each   sector is the same, as shown in Figure~\ref{fig:lattice} (a), up to a reflection on a line through the origin, since the interior of the sector contains no equilibria and the invariant circle is globally attracting by Theorem~\ref{IST}. 
Hence the part of the cell complex corresponding to the sector is always the same: two 1-dimensional cells, six separatrices giving rise to 0-dimensional cells with 1-dimensional fibre, five equilibria yielding  0-dimensional cells with 0-dimensional fibre, with the order shown in Figure~\ref{fig:lattice} (b).

The global cell complex is a concatenation of those obtained from the sectors, depending on the stability  within the invariant circle of the points denoted $p_1$ and $q_1$ in Figure~\ref{fig:lattice} (a).
In order to construct it, we start with the sector determined by $\theta_1$ and $\theta_2$.
The point $q_1$ is an attractor if and only if it determines a $(1-)$ in $\sigma(\GG)$.
Then the dynamics, and hence the cell complex, in the sector coming after this one is a reflection of that of Figure~\ref{fig:lattice} on the line containing the ray from the origin to $q_1$.
The other possibility is that $q_1$ is a saddle-node with symbol $(2+)$ in   $\sigma(\GG)$, and hence the sector coming after and its cell complex are copies of the first sector and its cell complex.

From the reasoning above it is clear that for $\GG_1=\LL(Q_1)$,  $\GG_2=\LL(Q_2)$, we have $\sigma(\GG_1)\equiv\sigma(\GG_2)$ if and only if they correspond to dynamics on ${\mathbb D}$ with isomorphic  cell complexes.
From \cite[Theorem 2']{NewmannOBrien}, two continuous flows on the plane with finitely many separatrices are topologically equivalent if and only if they have isomorphic cell complexes.
 It follows that $Q_1\sim Q_2$ if and only if $\sigma(\GG_1)\equiv\sigma(\GG_2)$.
\end{proof}

Thus, the global dynamics of  \eqref{eqR2contracting} for $Q\in \contr^{2p+1}$ is completely determined by the dynamics on the invariant circle, or equivalently, by the dynamics on the circle at infinity of the Poincar\'e disc.
When $Q$ is contracting the  dynamics of  \eqref{eqR2contracting}  only depends on  the polynomial $\LL Q$, in sharp contrast with the general (not contracting) case where the dynamics also depends on $\MM Q$, as described in \cite{ACL2021}.

The invariant may now be used to  decompose  $\contr^{2p+1}$ under $\sim$ into the following sets:
\begin{itemize}
\item[]
$\Sigma_0$ is the set of $Q\in\contr^{2p+1}$ such that $\sigma(\LL (Q))$ does not contain the symbols $(2\pm)$ and $\sigma(\LL (Q))\ne\infty$;
\item[]
$\Sigma_j$ for $j=1,2,\ldots, p+1$ is the set of $Q\in\contr^{2p+1}$ such that $\sigma(\LL (Q))$  contains exactly $j$ occurrences of the symbols $(2\pm)$;
\item[]
$\Sigma_{p+2}$ is the set of $Q\in\contr^{2p+1}$ such that $\sigma(\LL (Q))=\infty$.
\end{itemize}

The next result describes the geometry of these sets. 
In particular, it follows that generically $Q\in\Sigma_0$.

\begin{theorem}\label{teo:classifica}
The sets $\Sigma_j\subset\contr^{2p+1}$ satisfy:
\begin{enumerate}
\renewcommand{\theenumi}{(\alph{enumi})}
\renewcommand{\labelenumi}{{\theenumi}}
\item\label{item:sigma0}
$\Sigma_0$ is the union of an open and dense subset of $\contr^{2p+1}$ with a set of codimension 2 in $\contr^{2p+1}$;
\item\label{item:sigmaj}
each $\Sigma_j$, $j=1,2,\ldots, p+1$  is the union of  a subset of codimension $j$ of $\contr^{2p+1}$ and
 a set of codimension $2j+2$ in $\contr^{2p+1}$;
\item\label{item:sigmap+2}
$\Sigma_{p+2}=\contr^{2p+1}\cap\ker\LL$ and has codimension $2p+3$ in $\contr^{2p+1}$.
\end{enumerate}
\end{theorem}

\begin{proof}
The main argument in the proof is that for $A\subset  P^{2p+2}(\plano,\RR)$ we have  that $\cod \LL^{-1}(A)\cap \contr^{2p+1}=\cod(A)$.
This is true because $\contr^{2p+1}$ is an open subset of $ P^{2p+1}(\plano,\plano)$ and $\LL(\contr^{2p+1})= P^{2p+2}(\plano,\RR)$ by Theorem~\ref{teo:contract}.

The set $\OO_0$ of  polynomials that only have simple roots in $\projective$ is open and dense in $ P^{2p+2}(\plano,\RR)$.
Since $\LL$ is a continuous and open map, therefore $\LL^{-1}(\OO_0)\subset\Sigma_0$ is open and dense in $\contr^{2p+1}$.
The complement $\Sigma_0\backslash\LL^{-1}( \OO_0)$ consists of those $Q$ such that $\LL Q$ has at least one root of multiplicity at least 3 in $\projective$, and this latter set  is the union of sets of codimension $\ge 2$. 
This establishes \ref{item:sigma0}.

Similarly, the set $\OO_j$, $j=1,2,\ldots, p+1$ of polynomials with simple roots in $\projective$, except for exactly $j$ roots of multiplicity 2 satisfies $\cod\OO_j=j$ in $ P^{2p+2}(\plano,\RR)$ and 
$\LL^{-1}(\OO_j)\subset\Sigma_j$.
The complement $\Sigma_j\backslash\LL^{-1}( \OO_j)$ consists of those $Q$ such that either one of the roots of  $\LL Q$  in $\projective$ that corresponds to a symbol $(1\pm)$  has  multiplicity at least 3, or one of the roots corresponding to a symbol $(2\pm)$ has multiplicity at least 4, establishing \ref{item:sigmaj}.

Finally, $\Sigma_{p+2}=\ker\LL\cap\contr^{2p+1}$ and hence $\cod\Sigma_{p+2}=\dim P^{2p+2}(\plano,\RR)=2p+3$, hence \ref{item:sigmap+2} holds.
\end{proof}

The partition  $\contr^{2p+1}=\bigcup_{j=0}^{p+2}\Sigma_j$  is not a stratification of $\contr^{2p+1}$.
For instance, polynomials with  $\sigma(\GG)=(2+)(2+)$ and two different roots of multiplicity 2 may accumulate on a polynomial with a single root of multiplicity 4 for which $\sigma=(2+)$.
Therefore, the closure of $\Sigma_2$, a set of codimension 2, contains points of $\Sigma_1$ that has lower codimension.

\section{A class of examples --- definite nonlinearities}\label{sec:exDefinite}
We consider the family of planar vector fields given in \cite{GLS1987}
\begin{equation} \label{eq: ex sotomayor}
X(v)=A v+\va(v)Bv
\end{equation}
where $A=\begin{pmatrix}\lambda & 0 \\
                         0 & \lambda \\
        \end{pmatrix}$, $\lambda >0$, $B=\begin{pmatrix}
                                a & b \\
                                c & d \\
                              \end{pmatrix}$ is a $2\times 2$ matrix and $\va:\RR \to \RR$ is a homogeneous polynomial of even degree $2n\geq 2$.

The polynomial $\va $ is said to be positive (negative) definite if $\va(v)>0$ ($\va(v)<0$) for all $v\neq (0,0)$. Hence, the polynomial $Q(x,y)=\va(x,y)B \begin{pmatrix}
                                        x \\
                                        y \\
                                      \end{pmatrix}
$ is contracting provided by $\va(v)$ is positive (negative) definite and $B$ is a negative (positive) definite matrix in the sense that $\xy B \xy^T$ is a negative (positive) definite binary form.
In that case we say that $\va $ and $B$ are {\em of opposite sign}.

\begin{figure}
  \begin{center}
  \parbox{33mm}{
    \begin{center}
  \includegraphics[width=30mm]{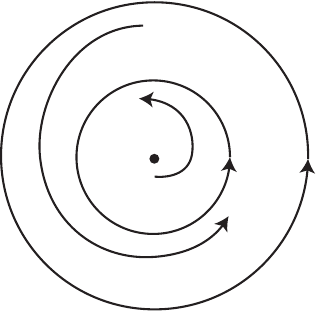}\\
\ref{Soto1} $\psi=\alpha(x^2+y^2)$
 \end{center}
  }
  \qquad
  \parbox{31mm}{
    \begin{center}
  \includegraphics[width=30mm]{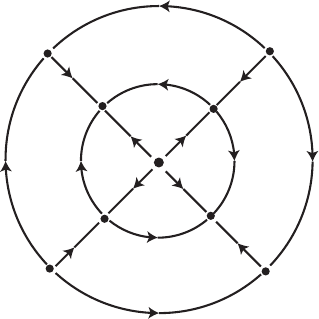}\\
\ref{Soto2} $\psi=x^2-y^2$
 \end{center}
  }
  \qquad
  \parbox{31mm}{
    \begin{center}
  \includegraphics[width=30mm]{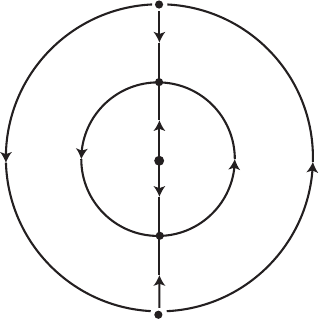}\\
\ref{Soto3} $\psi=\alpha x^2$
 \end{center}
  }
  \qquad
  \parbox{31mm}{
    \begin{center}
  \includegraphics[width=30mm]{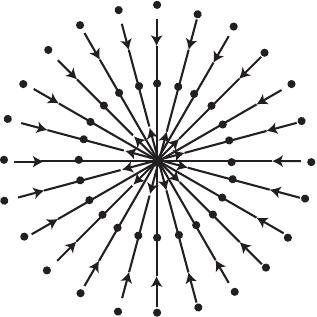}\\
\ref{Soto4} $\psi=0$
 \end{center}
  } 
  \end{center}
  \caption{Phase-portraits for Proposition \ref{prop: ex sotomayor} where $\alpha=\pm 1$.}\label{fig: ex sotomayor}
\end{figure}

\begin{proposition} \label{prop: ex sotomayor}Suppose $\va $ and $B$ in \eqref{eq: ex sotomayor} are of opposite sign. Then, the qualitative phase-portrait of \eqref{eq: ex sotomayor} (up to orientation of the orbits) is of one of types given in Figure~\ref{fig: ex sotomayor}.
\end{proposition}

\begin{proof}
Since $\LL Q(x,y)=\va(x,y)(cx^2+(d-a)xy-by^2)$ and $\va $ is definite, the phase-portrait at infinity is given by the second order binary form $\psi(x,y)=cx^2+(d-a)xy-by^2$ (up the orientation of the orbits). According to \cite[Theorem 1.3]{AnnaL1990} and the proof of Theorem \ref{teo:cubicNforms} below, the following vector fields give all possible dynamics of the vector field \eqref{eq: ex sotomayor} (up to orientation of the orbits):

\begin{enumerate}
\renewcommand{\theenumi}{(\Roman{enumi})}
\renewcommand{\labelenumi}{{\theenumi}}
  \item\label{Soto1}
  $\left\{\begin{array}{lcl}
\dot{x} &=& \lambda x + \va(x,y)(ax-\alpha y)\\
\dot{y} &=& \lambda y + \va(x,y)(\alpha x+a y), \quad \alpha= \pm 1
\end{array}  \right.$\\

  \item\label{Soto2}
  $\left\{\begin{array}{lcl}
\dot{x} &=& \lambda x + \va(x,y)(ax+ y)\\
\dot{y} &=& \lambda y + \va(x,y)(x+a y), \quad a^2-1>0
\end{array}  \right.$ \\

  \item\label{Soto3}
  $\left\{\begin{array}{lcl}
\dot{x} &=& \lambda x + a\va(x,y)x\\
\dot{y} &=& \lambda y + \va(x,y)(\alpha x+a y), \quad \alpha= \pm 1, \quad 
4a^2-1>0
%a\neq 0
\end{array}  \right.$ \\

  \item\label{Soto4}
  $\left\{\begin{array}{lcl}
\dot{x} &=& \lambda x + a \va(x,y)x\\
\dot{y} &=& \lambda y + a \va(x,y)y, \quad a\neq 0
\end{array}  \right.$ \\
\end{enumerate}

If $a<0$ ($a>0$) in the normal forms above, then $B$ is negative (positive) definite. So, if $\va(x,y)a<0$ the vector field $Q(x,y)=\va\xy B \xy^T$ is contracting and by Theorem \ref{IST} there exists a globally attracting circle. 
The dynamics on the circle is given by $\dot{\theta}=\LL Q(\cos{\theta},\sin{\theta})=g(\theta)$ and coincides with the dynamics on the circle at infinity.
The expressions for $g(\theta)$ for \ref{Soto1}--\ref{Soto4}  are given in Table~\ref{table:Soto}, hence  the phase-portraits are those  in Figure~\ref{fig: ex sotomayor}.
\end{proof}

\begin{table}
\begin{tabular}{l|cccc}
case	&	\ref{Soto1}	&	\ref{Soto2}	&	\ref{Soto3}	&	\ref{Soto4}	\\
\hline
$g(\theta)$	&	$\alpha \va(\cos{\theta},\sin{\theta})$	&	$\cos{(2\theta)}\va(\cos{\theta},\sin{\theta})$	&	$\alpha \cos^2{\theta}\va(\cos{\theta},\sin{\theta})$	&	0	\\
&			&		&		&		\\
$\sigma(\LL Q)$&$\varnothing$&(1+)(1-)&(2+)&$\infty$\\
\multicolumn{2}{}			&		&		&		\\
\end{tabular}
\caption{Expressions of $g(\theta)=\LL Q(\cos{\theta},\sin{\theta})$ and symbol sequences for the normal forms in Proposition~\ref{prop: ex sotomayor}  $\sigma(\LL Q)$ where $\alpha=\pm 1$.}\label{table:Soto}
\end{table}

Observe that the family \eqref{eq: ex sotomayor} realizes all possibilities given by Proposition~\ref{prop:global_dynamics pi}. In case $\lambda <0$, the polynomial $\va $ and the matrix $B$ must be of the same sign for Proposition \ref{prop: ex sotomayor} to hold.

\section{Another class of examples --- cubic nonlinearities}\label{sec:cubic}
We can now describe the phase diagrams for star nodes in the plane with contracting homogeneous cubic nonlinearity.
First we note that Proposition~\ref{prop:ContractGeneral} takes a particularly simple form
stated in the next result:

\begin{corollary}[of Proposition~\ref{prop:ContractGeneral}]\label{lemmaContract:n=3}
A homogeneous polynomial vector field $Q$ of degree $3$  in $\plano$ is contracting if writing $Q$ in the form \eqref{formNonHomogeneous} the following conditions hold:
\begin{enumerate}
\renewcommand{\theenumi}{(\roman{enumi})}
\renewcommand{\labelenumi}{{\theenumi}}
\item\label{item:deg3contract1}
either $p_1(1,0)$ and $p_1(0,1)<0$ or $p_2(1,0)$ and $p_2(0,1)<0$, and
\item\label{item:deg3contract2}
$4p_1(1,0)p_2(1,0)>(p_3(1,0)+p_4(1,0))^2$ 
\item\label{item:deg3contract3} 
$4p_1(0,1)p_2(0,1)>(p_3(0,1)+p_4(0,1))^2$.
\end{enumerate}
\end{corollary}
\begin{proof}
Writing $p_j(x^2,y^2)=a_{j0}x^2+a_{j1}y^2$ for $j=1,\ldots,4$, the result follows from Proposition~\ref{prop:ContractGeneral}, since for all $(x,y)\neq (0,0)$:
\begin{itemize}
\item[\ref{item:deg3contract1}]
  implies that either $p_1(x^2,y^2)=a_{10}x^2+a_{11}y^2<0$ or $p_2(x^2,y^2)=a_{20}x^2+a_{21}y^2<0$). It has been observed before that together with the second condition (here ensured by \ref{item:deg3contract2} and \ref{item:deg3contract3}), it is enough to verify one of the inequalities.
\item[\ref{item:deg3contract2}] and \ref{item:deg3contract3}
 imply that $4p_1(x^2,y^2)p_2(x^2,y^2)-(p_3(x^2,y^2)+p_4(x^2,y^2))^2=(4a_{10}a_{20}-(a_{30}+a_{40})^2)x^4
      +(4a_{10}a_{21}+4a_{11}a_{20}-2(a_{30}+a_{40})(a_{31}+a_{41}))x^2y^2+(4a_{11}a_{21}-(a_{31}+a_{41})^2)y^4>0$ since all the coefficients are positive.
      Note that $4a_{10}a_{20}>(a_{30}+a_{40})^2$ and $4a_{11}a_{21}>(a_{31}+a_{41})^2$ imply that $(a_{30}+a_{40})(a_{31}+a_{41}) < 4 \sqrt{a_{10}a_{20}a_{11}a_{21}}$ which is smaller than $2(a_{10}a_{21}+a_{11}a_{20})$.
\end{itemize}
\end{proof}

%%
% novo
%%

\begin{table}[htt]
\begin{center}
\renewcommand{\thecubico}{(\Roman{cubico})}
\begin{tabular}{cll}
  \refstepcounter{cubico}	\thecubico\label{Contract1}
  &$\left\{\begin{array}{lcl}
\dot{x} &=& \lambda x +3\mu x(x^2+y^2)-y^3 \\
\dot{y} &=& \lambda y  +3\mu y(x^2+y^2)+x(x^2+6\mu y^2)
\end{array}  \right.
$&$\begin{array}{l}
\mu<-\frac{1}{3}\\
{\mathcal G}\xy=x^4+6\mu x^2y^2+y^4
\end{array}$\\
\\
  \refstepcounter{cubico}	\thecubico\label{Contract2}
  &$\left\{\begin{array}{lcl}
\dot{x} &=& \lambda x -Kx(x^2+y^2)-\alpha y^3 \\
\dot{y} &=& \lambda y  -Ky(x^2+y^2)-\alpha x(x^2+6\mu y^2)
\end{array}  \right.
$&$\begin{array}{l}
\alpha=\pm 1, \mu>-\frac{1}{3}, \mu \neq \frac{1}{3}\\
{\mathcal G}\xy=\alpha (x^4+6\mu x^2y^2+y^4)\\
K=\max\{(3\mu)^2,1/2\}
\end{array}$\\
\\
  \refstepcounter{cubico}	\thecubico\label{Contract3}
  &$\left\{\begin{array}{lcl}
\dot{x} &=& \lambda x -Kx(x^2+y^2) +y^3\\
\dot{y} &=& \lambda y  -Ky(x^2+y^2)+x(x^2+6\mu y^2)
\end{array}  \right.$
& $\begin{array}{l}
{\mathcal G}\xy=x^4+6\mu x^2y^2-y^4\\
K=\max\{(3\mu)^2,1/2\}
\end{array}$ \\
\\
  \refstepcounter{cubico}	\thecubico\label{Contract4}
  &$\left\{\begin{array}{lcl}
\dot{x} &=& \lambda x -4x(x^2+y^2)-\alpha y(6x^2+y^2) \\
\dot{y} &=& \lambda y  -4y(x^2+y^2)
\end{array}  \right.
$&$\begin{array}{l}
 \alpha=\pm 1\\
{\mathcal G}\xy=\alpha(6x^2y^2+y^4)
\end{array}$\\
\\
  \refstepcounter{cubico}	\thecubico\label{Contract5}
  &$\left\{\begin{array}{lcl}
\dot{x} &=& \lambda x -x(x^2+y^2)-\alpha y(+x^2-y^2)\\
\dot{y} &=& \lambda y -y(x^2+y^2)
\end{array}  \right.
$&$\begin{array}{l}
 \alpha=\pm 1\\
{\mathcal G}\xy=\alpha(6x^2y^2-y^4)
\end{array}$\\
\\
%  \refstepcounter{cubico}	\thecubico\label{Contract6}
%  &$\left\{\begin{array}{lcl}
%\dot{x} &=& \lambda x -x(x^2+y^2)-\alpha y^3\\
%\dot{y} &=& \lambda y -y(x^2+y^2)+\alpha x(x^2+2y^2)
%\end{array}  \right.
%$&$\begin{array}{l}
% \alpha=\pm 1\\
%{\mathcal G}\xy=\alpha(x^2+y^2)^2
%\end{array}$\\
%\\
  \refstepcounter{cubico}	\thecubico\label{Contract7}
  &$\left\{\begin{array}{lcl}
\dot{x} &=& \lambda x -x(x^2+y^2)\\
\dot{y} &=& \lambda y -y(x^2+y^2)+6xy^2
\end{array}  \right.$
& ${\mathcal G}\xy=6x^2y^2$ \\
\\
%  \refstepcounter{cubico}	\thecubico\label{Contract8}
%  &$\left\{\begin{array}{lcl}
%\dot{x} &=& \lambda x -2x(x^2+y^2)\\
%\dot{y} &=& \lambda y -2y(x^2+y^2)+4yx^2
%\end{array}  \right.$
%& ${\mathcal G}\xy=4yx^3$ \\
%\\
%  \refstepcounter{cubico}	\thecubico\label{Contract9}
%  &$\left\{\begin{array}{lcl}
%\dot{x} &=& \lambda x -x(x^2+y^2)\\
%\dot{y} &=& \lambda y -y(x^2+y^2)+\alpha x^3
%\end{array}  \right.
%$&$\begin{array}{l}
%\alpha= \pm 1\\
%{\mathcal G}\xy=\alpha x^4
%\end{array}$\\
%\\
  \refstepcounter{cubico}	\thecubico\label{Contract10}
  &$\left\{\begin{array}{lcl}
\dot{x} &=& \lambda x -x(x^2+y^2)\\
\dot{y} &=& \lambda y -y(x^2+y^2)
\end{array}  \right.$
& ${\mathcal G}\xy=0$ \\
\setcounter{lixo}{\value{cubico}}
\end{tabular}
\end{center}
\caption{Normal forms for contracting cubic nonlinearities in Theorem~\ref{teo:cubicNforms} with $\GG\xy=\LL Q\xy$.}\label{table:Cubico}
\end{table}

\begin{figure}
\begin{center}
\parbox{30mm}{
\begin{center}\includegraphics[width=30mm]{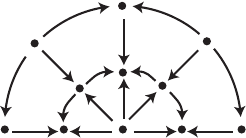}\\
\ref{Contract1}
\end{center}}
\qquad
\parbox{30mm}{
\begin{center}\includegraphics[width=30mm]{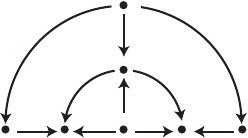}\\
\ref{Contract3} 
%and \ref{Contract8}
\end{center}}
\qquad
\parbox{30mm}{
\begin{center}\includegraphics[width=30mm]{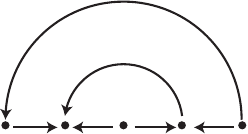}\\
\ref{Contract4} 
%and \ref{Contract9} 
\end{center}}
\\
\parbox{30mm}{
\begin{center}\includegraphics[width=30mm]{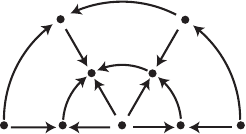}\\
\ref{Contract5}
\end{center}}
\qquad
\parbox{30mm}{
\begin{center}\includegraphics[width=30mm]{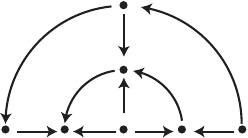}\\
\ref{Contract7}
\end{center}}
\qquad
\parbox{30mm}{
\begin{center}\includegraphics[width=30mm]{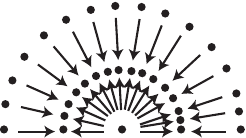}\\
\ref{Contract10}
\end{center}}
\\
\parbox{30mm}{
\begin{center}\includegraphics[width=30mm]{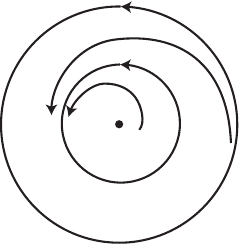}\\
\ref{Contract2} 
%and \ref{Contract6}
\end{center}}
\end{center}
\caption{Qualitative portraits on the Poincar\'e disc for  \eqref{eqR2contracting} with a contracting cubic nonlinearity.
 On the six top figures only half the disc is shown, the other half is obtained by rotation of $\pi$ around the origin.
Numbering corresponds to normal forms in Theorem~\ref{teo:cubicNforms}.
In Case~\ref{Contract10} both the circle at infinity and the invariant sphere are continua of equilibria. }\label{fig: teo cubic}
\end{figure}

\begin{theorem} \label{teo:cubicNforms} 
Let $\lambda>0$ and $Q$ be a contracting homogeneous cubic vector field.
Then \eqref{eqR2contracting} is equivalent to one of the 7 normal forms
in Table~\ref{table:Cubico}.
% below
The qualitative phase-portrait of \eqref{eqR2contracting} is of one of types shown in Figure~\ref{fig: teo cubic} and symbol sequences are given in Table~\ref{table:cubicNforms}.
\end{theorem}

\begin{table}[hhh]
\begin{center}
\setcounter{cubico}{\value{lixo}}
\renewcommand{\thecubico}{(\Roman{cubico})}
\begin{tabular}{cll}
  \refstepcounter{cubico}	\thecubico\label{Contract6}
  &$\left\{\begin{array}{lcl}
\dot{x} &=& \lambda x -x(x^2+y^2)-\alpha y^3\\
\dot{y} &=& \lambda y -y(x^2+y^2)+\alpha x(x^2+2y^2)
\end{array}  \right.
$&$\begin{array}{l}
 \alpha=\pm 1\\
{\mathcal G}\xy=\alpha(x^2+y^2)^2
\end{array}$\\
\\
  \refstepcounter{cubico}	\thecubico\label{Contract8}
  &$\left\{\begin{array}{lcl}
\dot{x} &=& \lambda x -2x(x^2+y^2)\\
\dot{y} &=& \lambda y -2y(x^2+y^2)+4yx^2
\end{array}  \right.$
& ${\mathcal G}\xy=4yx^3$ \\
\\
  \refstepcounter{cubico}	\thecubico\label{Contract9}
  &$\left\{\begin{array}{lcl}
\dot{x} &=& \lambda x -x(x^2+y^2)\\
\dot{y} &=& \lambda y -y(x^2+y^2)+\alpha x^3
\end{array}  \right.
$&$\begin{array}{l}
\alpha= \pm 1\\
{\mathcal G}\xy=\alpha x^4
\end{array}$\\
\\
\end{tabular}
\end{center}
\caption{Normal forms for contracting cubic nonlinearities in the list of \cite[Theorem 2.6]{AnnaL1990} that do not appear in Theorem~\ref{teo:cubicNforms} with $\GG\xy=\LL Q\xy$.}\label{table:CubicoNao}
\end{table}

\begin{table}[hht]
\begin{tabular}{cclllc}
normal	&	equilibria	&	type of	&	angular	&	symbol	&	$\cod\Sigma_j$	 \\
form	&	at infinity	&	roots	&	stability	&	sequence	&		 \\ \hline
\ref{Contract1}	&	8	&	simple	&	hyperbolic	&	$(1+),(1-),(1+),(1-)$	&	0	 \\ \hline
\ref{Contract2}	&	0	&	none	&	-	&	$\varnothing$	&	0	 \\ \hline
\ref{Contract3}	&	4	&	simple	&	hyperbolic	&	$(1+),(1-)$	&	0	 \\ \hline
\ref{Contract4}	&	2	&	double	&	saddle-nodes	&	$(2+)$	&	1	 \\ \hline
\ref{Contract5}	&	6	&	4 simple	&	4 hyperbolic	&	$(2+),(1-),(1+)$	&	1	 \\ 
	&		&	2 double	&	2 saddle-nodes	&	or $(2-),(1+),(1-)$	&		 \\ \hline
\ref{Contract7}	&	4	&	double	&	 saddle-nodes	&	$(2+),(2+)$	&	2	\\
	&	(0, $\pi/2$, $\pi$, $3\pi/2$)	&		&		&		&		 \\ \hline
\ref{Contract10}	&	$\infty$	&	all	&	-	&	$\infty$	&	3	 \\ \hline \hline
\ref{Contract6}	&	0	&	none	&	-	&	$\varnothing$	&	0	 \\ \hline	
\ref{Contract8}	&	4	&	2 simple	&	2 hyperbolic	&	$(1+),(1-)$	&	0	 \\
	&	(0, $\pi/2$, $\pi$, $3\pi/2$)	&	2 triple	&	2 hyperbolic-like	&		&		 \\ \hline
\ref{Contract9}	&	2	&	quadruple	&	saddle-nodes	&	$(2+)$	&	1	\\
	&	($\pi/2$, $3\pi/2$)	&		&		&		&		 \\ \hline	\\
	\end{tabular}
\caption{Number and angular stability of finite and infinite equilibria for normal forms in Theorem~\ref{teo:cubicNforms}. Hyperbolic-like are weak non-hyperbolic attractors or repellors. Symbol sequences refer to the coding of Section~\ref{subsec:classification}  and the sign $(2+)$
may be replaced by $(2-)$ depending on the value of $\alpha=\pm 1$. The subset $\Sigma_j$ is that of Theorem~\ref{teo:classifica}.}\label{table:cubicNforms}
\end{table}

\begin{proof}
Normal forms for binary forms of degree 4, up to a linear change of coordinates, are given in  \cite[Theorem 2.6]{AnnaL1990}.
For each binary form $\GG$ on this list,  Theorem~\ref{teo:contract} ensures that there is a vector field  \eqref{eqR2contracting} with contracting nonlinearity $Q$ such that $\GG\xy=\LL Q\xy$.
Since the dynamics of  \eqref{eqR2contracting} is totally determined by $g(\theta)=\GG(\cos\theta,\sin\theta)$ and since by Lemma~\ref{lem:linearChange} a linear change of coordinates in \eqref{eqR2contracting}  corresponds to a linear change of coordinates in $\GG$,
this gives a list of all possible dynamical behaviour.

The list of  \cite[Theorem 2.6]{AnnaL1990} contains ten normal forms, three of which do not appear in our list because they yield dynamics that is globally equivalent to one of the forms in Table~\ref{table:Cubico}. 
They are listed in Table~\ref{table:CubicoNao}.

The cubic nonlinearities  in both lists were obtained following the construction in the proof of  Theorem~\ref{teo:contract}.
The constant $K$ such that $Q$ is contracting was obtained from Corollary~\ref{lemmaContract:n=3} as follows: in all cases, except for \ref{Contract8}, the binary form $\GG$ is written as $\GG\xy=x^2b_1(x^2,y^2)+y^2b_3(x^2,y^2)$.
This yields, in the notation of Proposition~\ref{prop:homogeneous_polynomial}, the choices $p_1(x^2,y^2)=p_2(x^2,y^2)=-K(x^2+y^2)$, hence $p_1(1,0)=p_2(1,0)=p_1(0,1)=p_2(0,1)=-K<0$.
Conditions \ref{item:deg3contract2} and \ref{item:deg3contract3} of the corollary become
$$
4K^2>(p_3(1,0)+p_4(1,0))^2
\quad\mbox{and}\quad
4K^2>(p_3(0,1)+p_4(0,1))^2 .
$$
These expressions are evaluated in Table~\ref{table:calculaK}.
For the remaining case \ref{Contract8} we have $p_1(x^2,y^2)=-K(x^2+y^2)$ and $p_2(x^2,y^2)=-K(x^2+y^2)+4x^2$ with $p_3(x^2,y^2)=p_4(x^2,y^2)=0$. 
Conditions \ref{item:deg3contract2} and \ref{item:deg3contract3} of the corollary are then 
$-4K(4-K)>0$ and $4K^2>0$, satisfied by any $K\in (0,4)$, for instance, $K=2$.

Since $\lambda>0$ and the nonlinearities in Systems \ref{Contract1}--\ref{Contract10} are contracting,
it follows by Theorem \ref{IST} that there exists a globally attracting circle. 
The dynamics on the circle is given by  $\dot \theta=g(\theta)$, where $g(\theta)=\mathcal{G}(\cos{\theta},\sin{\theta})$ 
and coincides with the dynamics on the circle at infinity.
The number of solutions of $g(\theta)=0$, their type and stability are given in Table~\ref{table:cubicNforms}.

From Table~\ref{table:cubicNforms} it follows that the three normal forms in Table~\ref{table:CubicoNao} that share the same symbol sequence.
These are: 
$$
 \ref{Contract6}\sim \ref{Contract2}
 \qquad
 \ref{Contract8}\sim \ref{Contract3}
 \quad\mbox{and}\quad
 \ref{Contract9}\sim \ref{Contract4} .
 $$
Thus, these normal forms have globally equivalent dynamics.
\end{proof}

\begin{table}[hhh]
\begin{tabular}{ccccc}
normal	&	modal	&	$\left(p_3(1,0)+p_4(1,0)\right)^2$	&	$\left(p_3(0,1)+p_4(0,1)\right)^2$	&	$K^2$	 \\
form	&	parameter	&		&		&		 \\ \hline
\ref{Contract1}	&	$\mu<-1/3$	&	1	&	$(-1+6\mu)^2>4(3\mu)^2>4$	&	$(3\mu)^2$	 \\ \hline
\ref{Contract2}	&	$1/3\ne\mu>-1/3$	&	$\alpha^2=1$	&	$(1+6\mu)^2>4(3\mu)^2$	&	$\max\{(3\mu)^2,1/2\}$	 \\ 
	&	$\alpha=\pm 1$	&		&		&		 \\ \hline
\ref{Contract3}	&	$\mu\in\RR$	&	1	&	$(-1+6\mu)^2>4(3\mu)^2$	&	$\max\{(3\mu)^2,1/2\}$	 \\ \hline
\ref{Contract4}	&	$\alpha=\pm 1$	&	$(-6\alpha)^2$	&	$(-\alpha)^2=1$	&	16	 \\ \hline
\ref{Contract5}	&	$\alpha=\pm 1$	&	$(-\alpha)^2=1$	&	$\alpha^2=1$	&	1	 \\ \hline
\ref{Contract7}	&		&	0	&	1	&	1	 \\ \hline
\ref{Contract10}	&		&	0	&	0	&	1	 \\ \hline

\ref{Contract6}	&	$\alpha=\pm 1$	&	$\alpha^2=1$	&	$(-\alpha+2\alpha)^2=1$	&	1	 \\ \hline
\ref{Contract9}	&	$\alpha=\pm 1$	&	$\alpha^2=1$	&	0	&	1	 \\ \hline	\\	
	
	\end{tabular}
\caption{Calculation of a value of $K$ for which $Q$ in the normal forms in Theorem~\ref{teo:cubicNforms} is contracting from Corollary~\ref{lemmaContract:n=3}. The case \ref{Contract8} is different and $K$ is computed in the text.}\label{table:calculaK}
\end{table}

\subsection{Cubic  $\Ztwo\oplus\Ztwo$ nonlinearities}\label{sub:Z2}

As in \ref{subsec:Z2dyn} above, we may say more in the symmetric case,
writing
\begin{equation}\label{eq:QZ2}
Q\xy=\left( -x(a_{10}x^2+a_{11}y^2),-y(a_{20}x^2+a_{21}y^2)\right) 
\end{equation}
we start by finding when $Q$ is contracting.

\begin{theorem}\label{ThContractionDeg3}
The polynomial $Q$ of the form \eqref{eq:QZ2} is contracting if and only if $a_{10} >0$, $a_{21}>0$ and one of the following conditions holds:
\begin{enumerate}
\renewcommand{\theenumi}{(\alph{enumi})}
\renewcommand{\labelenumi}{{\theenumi}}
\item\label{CondbMaisc} 
$a_{11} +a_{20} \ge 0$;
\item\label{CondDet} 
$a_{10} a_{21} -(a_{11} +a_{20} )^2/4>0$.
\end{enumerate}
\end{theorem}
\begin{proof}
In this case we have
$$
\MM Q\xy=-x^2(a_{10}x^2+a_{11}y^2)-y^2(a_{20}x^2+a_{21}y^2).
$$
{\sl Sufficiency:}
If $a_{10} >0$, $a_{21} >0$ and \ref{CondbMaisc} holds then clearly $\MM Q\xy<0$ for all $\xy\ne (0,0)$.
%The case when \ref{CondDet} holds 
The case  \ref{CondDet} follows from the proof of Proposition~\ref{prop:ContractGeneral}.

{\sl Necessity:}
The function $\MM Q\xy$  is a quadratic form $p(u,v)$ on $(u,v)=(x^2,y^2)$, represented by the symmetric matrix
$$
M=\begin{pmatrix}-a_{10} &-(a_{11} +a_{20} )/2\\-(a_{11} +a_{20} )/2&-a_{21} \end{pmatrix}.
$$
If $Q$ is contracting then $\MM Q\xy<0$ for all $\xy\ne (0,0)$.
In particular,  $\MM Q(x,0)=-a_{10}x^4<0$ for $x\ne 0$ and $\MM Q(0,y)=-a_{20}y^4<0$ for $y\ne 0$, hence $a_{10} >0$ and $a_{21} >0$ and $\tr M<0$.

The condition $\MM Q\xy<0$ for all $\xy\ne (0,0)$ implies that
$p(u,v)=-a_{10} u^2-(a_{11} +a_{20} )uv-a_{21} v^2<0$ for all $u\ge 0$, $v\ge 0$ with $(u,v)\ne (0,0)$.
Let $\mu_+\ge\mu_-$ be the eigenvalues of $M$.
Since $\tr M<0$ then $\mu_-<0$.
There are three possibilities:
\begin{enumerate}
\renewcommand{\theenumi}{(\roman{enumi})}
\renewcommand{\labelenumi}{{\theenumi}}
\item\label{menosMenos}
The quadratic form $p(u,v)$ is  negative definite, or equivalently, both  $\mu_+<0$ and $\mu_-<0$.
This implies $\det M>0$, hence \ref{CondDet} holds.
\item\label{menosZero}
The eigenvalues of $M$ satisfy  $\mu_-<0$ and $\mu_+=0$.
In suitable coordinates $(\tilde{u},\tilde{v})$, we have $p(\tilde{u},\tilde{v})=\mu_- \tilde{u}^2$, where $\tilde{u}$ is the coordinate in the direction of the eigenvector of $\mu_-$ and $\tilde{v}$ is the coordinate in the direction of the eigenvector of zero.
Thus, if $Q$ is contracting then the  eigenvector of zero does not lie in the first or third quadrants.
The eigenvectors $(u,v)$ of the zero eigenvalue satisfy $-a_{10} u-(a_{11} +a_{20} )v/2=0$, then they are scalar multiples of $(u,v)=(a_{11} +a_{20} ,-2a_{10} )$. This last vector is not in the first or the third quadrants if and only if $a_{11} +a_{20} >0$, as in \ref{CondbMaisc}.
\item\label{menosMais}
The eigenvalues of $M$  satisfy  $\mu_-<0$ and $\mu_+>0$.
In suitable coordinates  $(\tilde{u},\tilde{v})$, we have $p(\tilde{u},\tilde{v})=\mu_+ \tilde{u}^2+\mu_-\tilde{v}^2$, where $\tilde{u}$ is the coordinate the direction of the eigenvector of $\mu_+$ and $\tilde{v}$ is the coordinate in the direction of the eigenvector of $\mu_-$.
Therefore, if $Q$ is contracting, then the  eigenvector of $\mu_+$ does not lie in the (closure of) first nor in the third quadrant.

The characteristic polynomial of $M$ is
$$
p_M(\mu)=\mu^2+(a_{10} +a_{21} )\mu+a_{10} a_{21}-(a_{11} +a_{20} )^2/4
$$
and $2\mu_+=-(a_{10} +a_{21} )+\sqrt{\Delta}$ with $\Delta=(a_{10} -a_{21} )^2+(a_{11} +a_{20} )^2$.
The eigenvectors $(u,v)$ of $\mu_+$ satisfy $-a_{10} u-(a_{11} +a_{20} )v/2=\mu_+u$ or, equivalently,
$$
(a_{11} +a_{20} )v=\left(-2a_{10} -2\mu_+\right) u=
\left(a_{21} -a_{10} -\sqrt{\Delta}\right)u
$$
and are scalar multiples of $(u,v)=\left(a_{11} +a_{20} ,a_{21} -a_{10} -\sqrt{\Delta}\right)$.
Suppose $a_{11} +a_{20} <0$. Then we must have $a_{21} -a_{10} >\sqrt{\Delta}>0$ and this is equivalent to $(a_{21} -a_{10} )^2>\Delta=(a_{10} -a_{21} )^2+(a_{11} +a_{20} )^2$, or equivalently $0>(a_{11} +a_{20} )^2$, a contradiction.
Hence, $a_{11} +a_{20} \ge 0$.
\end{enumerate}
\end{proof}

From \eqref{eq:QZ2}  we get:
\begin{equation}\label{eq:LQZ2}
\LL Q\xy=xy\left(Ax^2-By^2\right)=xyq\xy
\qquad
A=a_{10} -a_{20} \quad\mbox{and}\quad B=a_{21}-a_{11} .
\end{equation}
The dynamics is completely determined by the values of $A$ and $B$, as the next result shows.

\begin{figure}[hht]
\begin{center}
    \includegraphics[width=40mm]{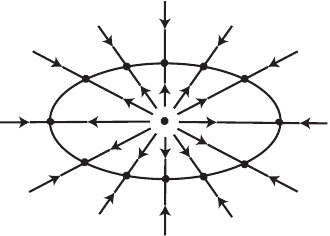}
    \end{center}
  \caption{Planar dynamics as in Proposition~\ref{prop:QZ2} \ref{item:Z21}. The invariant sphere is an ellipse of equilibria.}\label{grafradialsphere}
\end{figure}

\begin{proposition}\label{prop:QZ2}
If $Q$ is a contracting cubic  $\Ztwo\oplus\Ztwo$ equivariant vector field then the dynamics of  \eqref{eqR2contracting} in the invariant circle is the following:
\begin{enumerate}
\renewcommand{\theenumi}{(\Roman{enumi})}
\renewcommand{\labelenumi}{{\theenumi}}
\item\label{item:Z21}
If $A=B=0$ then $\sigma(\LL Q)=\infty$.
\item\label{item:Z22}
If $AB=0$ and $A+B\ne 0$ then $\sigma(\LL Q)=(1+)(1-)$ and one of the equilibria is not hyperbolic.
\item\label{item:Z23}
If $AB<0$ then $\sigma(\LL Q)=(1+)(1-)$.
\item\label{item:Z24}
If $AB>0$ then $\sigma(\LL Q)=(1+)(1-)(1+)(1-)$.
\end{enumerate}
Moreover, in cases \ref{item:Z23} and \ref{item:Z24} all equilibria  are hyperbolic.
\end{proposition}
\begin{proof}
First note that from \eqref{eq:LQZ2} there are always equilibria on the axes, at the 4 points where they cross the invariant circle.
Equilibria on the invariant circle are hyperbolic if and only if they are simple roots of $\LL Q$.

\ref{item:Z21} From Lemma~\ref{corContinuum}  the invariant circle is a continuum of equilibria if and only if $A=B=0$, establishing \ref{item:Z21}.
The invariant circle is the ellipse $a_{10} x^2+a_{11} y^2=\lambda$,  all the trajectories are contained in lines through the origin and go from the origin (or from infinity) to a point in the ellipse.
Indeed, $\dfrac{\dot y}{\dot x}=\dfrac{y}{x}$, hence $\dfrac{dy}{dx}=\dfrac{y}{x}$ and $y=Kx$, where $K$ is a real constant. See Figure~\ref{grafradialsphere}.

\ref{item:Z22}
If $A\ne 0$ and $B= 0$ then $q\xy=Ax^2$ so all the equilibria lie on the axes.  
The equilibria on the $x=0$ axis  not hyperbolic, since they are roots of multiplicity 3 of $\LL Q$.
The case $A= 0$ and $B\ne 0$ is analogous.

When both $A\ne 0$ and $B\ne 0$ then $q(x,0)\ne 0 \ne q(0,y)$.
Therefore the equilibria on the axes are hyperbolic. 
Other equilibria satisfy $Ax^2=By^2$.
There are two cases to consider.

\ref{item:Z23}
If $AB<0$ then $q\xy=0$ has no solutions so all the equilibria lie on the axes.

\ref{item:Z24}
If $AB>0$ then $q\xy=0$ has  solutions $y=\pm\sqrt{Ax^2/B}$, corresponding to one hyperbolic equilibrium on the interior of each one of the quadrants in the plane.
\end{proof}

\subsection*{Acknowledgements:} The authors are grateful to P. Gothen, R. Prohens and A. Teruel for  fruitful conversations.
% The first author was partially supported by the Spanish Research Project MINECO-18-MTM2017-87697-P. The last two authors were partially supported by CMUP (UID/MAT/00144/2019), which is funded by FCT (Portugal) with national (MCTES) and European structural funds (FEDER), under the partnership agreement PT2020.

\end{document}